%% file: Final_Platoon_TAC_2011Dec17.tex
\documentclass[journal]{ieeetran}

\input{commands}

\usepackage{amssymb,latexsym,color,amsmath,pifont,epsfig,mathtools,dsfont,multirow,booktabs}
\usepackage{graphicx}
\usepackage{subfig}      
\usepackage{setspace,cite}

\input umsa.fd
\input umsb.fd
\input ulasy.fd

\newcommand{\bbR}{\mathbb{R}}
\newcommand{\mrd}{\mathrm{d}}

\newcommand{\mre}{\mathrm{e}}

\newcommand{\non}{\nonumber}
\newcommand{\veps}{\varepsilon}
\newcommand{\eS}{{\cal S}}

\newcommand{\bal}{\begin{align}}
\newcommand{\eal}{\end{align}}
\newcommand{\ds}{\displaystyle}
\newcommand{\bfo}{\mathbf{1}}

\newcommand{\Pctr}{\Pi_{\rm ctr}}

\newcommand{\tc}{\textcolor}

\IEEEoverridecommandlockouts                              

\begin{document}

\title{\huge \bf Optimal Control of Vehicular Formations with
Nearest Neighbor Interactions}

\author{Fu Lin,~\IEEEmembership{Student Member, IEEE\/}, Makan Fardad, and Mihailo R. Jovanovi\'c,~\IEEEmembership{Member, IEEE}
\thanks{Financial support from the National Science Foundation under
CAREER Award CMMI-06-44793 and under Awards CMMI-09-27720 and
CMMI-09-27509 is gratefully acknowledged.}
\thanks{
F.\ Lin and M.\ R.\ Jovanovi\'c are with the Department of
Electrical and Computer Engineering, University of Minnesota,
Minneapolis, MN 55455. M.\ Fardad is with the Department of
Electrical Engineering and Computer Science, Syracuse University, NY
13244. E-mails: fu@umn.edu, makan@syr.edu, mihailo@umn.edu.}}

\maketitle

    \begin{abstract}
We consider the design of optimal localized feedback gains for one-dimensional formations in which vehicles only use information from their immediate neighbors. The control objective is to enhance coherence of the formation by making it behave like a rigid lattice. For the single-integrator model with symmetric gains, we establish convexity, implying that the globally optimal controller can be computed efficiently. We also identify a class of convex problems for double-integrators by restricting the controller to symmetric position and uniform diagonal velocity gains. To obtain the optimal non-symmetric gains for both the single- and the double-integrator models, we solve a parameterized family of optimal control problems ranging from an easily solvable problem to the problem of interest as the underlying parameter increases. When this parameter is kept small, we employ perturbation analysis to decouple the matrix equations that result from the optimality conditions, thereby rendering the unique optimal feedback gain. This solution is used to initialize a homotopy-based Newton's method to find the optimal localized gain. To investigate the performance of localized controllers, we examine how the coherence of large-scale stochastically forced formations scales with the number of vehicles. We establish several explicit scaling relationships and show that the best performance is achieved by a localized controller that is both non-symmetric and spatially-varying.
    \end{abstract}


\section{Introduction}

\subsection{Background}

\begin{figure}
    \begin{center}
    \includegraphics[width=0.4\textwidth]{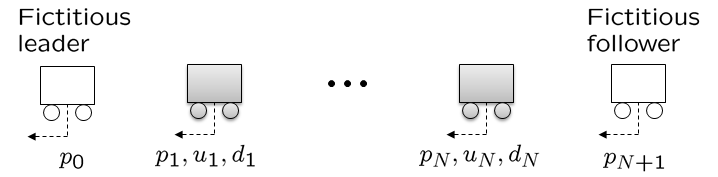}
    \caption{One-dimensional formation of vehicles.}
    \label{fig.1Dplatoon}
    \end{center}
    \end{figure}

The control of vehicular platoons has attracted considerable attention since the mid sixties~\cite{levath66,melkuo71,swahed96}. Recent technological advances in developing vehicles with communication and computation capabilities have spurred renewed interest in this area~\cite{jadlinmor03,seipanhed04,faxmur04,jovbam05,lafwilcauvee05,jovfowbamdan08,barmehhes09,midbra10,bamjovmitpatTAC12}. The simplest control objective for the one-dimensional (1D) formation shown in Fig.~\ref{fig.1Dplatoon} is to maintain a desired cruising velocity and to keep a pre-specified constant distance between neighboring vehicles. This problem is emblematic of a wide range of technologically relevant applications including the control of automated highways, unmanned aerial vehicles, swarms of robotic agents, and satellite constellations.

    \begin{table*}[ttt]
	\caption{\textup{Summary of asymptotic scalings with the number of vehicles $N$ for the optimal symmetric and non-symmetric position gains. The $N$-independent control penalty, $R = r_0 I$, in the quadratic performance objective leads to similar growth with $N$ of formation coherence and control energy (per vehicle). On the other hand, the $N$-dependent control penalty that provides bounded control energy yields less favorable coherence.}}
			\label{tab.main}
	\centering
	\begin{tabular}{|c||c|c|c|}
				\toprule
    \begin{tabular}{c}
	\textbf{Optimal}
    \\
    \textbf{position gains}
    \end{tabular}
    &
    \begin{tabular}{c}
    \textbf{Control penalty }
    \\[0.cm]
    \textbf{$R \, = \, r I$}
    \end{tabular}
	&
    \begin{tabular}{c}
    \textbf{Control energy}
    \\[0.cm]
    \textbf{(per vehicle)}
    \end{tabular}
    &
    \begin{tabular}{c}
    \textbf{Formation}
    \\
    \textbf{coherence}
    \end{tabular}
    \\
	\midrule
	\textbf{symmetric}
    &
    \multirow{2}{*}{$ \ba{rcl} & & \\[-0.35cm] r (N) & = & r_0 \\ & = & \mbox{\rm const.} \ea $}
    &
    $O(\sqrt{N})$
    &
    $O(\sqrt{N})$
    \\
    \cmidrule[0.55pt]{1-1} \cmidrule[0.55pt]{3-4}
    \textbf{non-symmetric}
    &
    &
    $O(\sqrt[4]{N})$
    &
    $O(\sqrt[4]{N})$
    \\
    \midrule
    \textbf{symmetric}
    &
    $r(N) \, \sim \, N$
    &
    $O(1)$
    &
    $O(N)$
    \\
    \midrule
    \textbf{non-symmetric}
    &
    $r(N) \, \sim \, \sqrt{N}$
    &
    $O(1)$
    &
    $O(\sqrt{N})$							
    \\
	\bottomrule	
	\end{tabular}
\end{table*}

Recent work in this area has focused on fundamental performance limitations of both centralized and decentralized controllers for large-scale formations~\cite{seipanhed04,jovbam05,jovfowbamdan08,barmehhes09,bamjovmitpatTAC12,midbra10}. For centralized linear quadratic optimal control formulations based on penalizing relative position errors it was shown in~\cite{jovbam05} that stabilizability and detectability deteriorate as formation size increases. In~\cite{jovfowbamdan08}, it was shown that merge and split maneuvers can exhibit poor convergence rates even upon inclusion of absolute position errors in cost functionals. In~\cite{seipanhed04}, it was shown that sensitivity of spacing errors to disturbances increases with the number of vehicles for formations with localized symmetric controllers that utilize relative position errors between neighboring vehicles. In~\cite{midbra10}, the analysis of~\cite{seipanhed04} was expanded to include heterogeneous vehicles, non-zero time headway, and limited communication range within the formation.

The motivation for the current study comes from two recent papers,~\cite{bamjovmitpatTAC12} and~\cite{barmehhes09}. In~\cite{bamjovmitpatTAC12}, fundamental performance limitations of {\em localized symmetric feedback\/} for spatially invariant consensus and formation problems were examined. It was shown that, in 1D, it is impossible to have coherent large formations that behave like rigid lattice. This was done by exhibiting linear scaling, with the number of vehicles, of the {\em formation-size-normalized\/} ${\cal H}_2$ norm from disturbances to an appropriately defined macroscopic performance measure. In 2D this measure increases logarithmically, and in 3D it remains bounded irrespective of the system size. These scalings were derived by imposing uniform bounds on control energy at each vehicle.

For formations on a one-dimensional lattice, it was shown in~\cite{barmehhes09} that the decay rate (with the number of vehicles) of the least damped mode of the closed-loop system can be improved by introducing a small amount of `mistuning' to the spatially uniform symmetric feedback gains. A large formation was modeled as a diffusive PDE, and an optimal small-in-norm perturbation profile that destroys the spatial symmetry and renders the system more stable was designed. Numerical computations were also used to demonstrate that the spatially-varying feedback gains have beneficial influence on the closed-loop ${\cal H}_\infty$ norm. The PDE approaches have also been found useful in the deployment of multi-agents~\cite{frikrs10a,frikrs10b} and in coordination algorithms~\cite{sarsep09}.

Even though traditional optimal control does not facilitate incorporation of structural constraints and leads to centralized architectures, the optimal feedback gain matrix for both spatially invariant systems~\cite{bampagdah02} and systems on graphs~\cite{motjad08} have {\em off-diagonal decay\/}. Several recent efforts have focused on identification of classes of convex distributed control problems. For spatially invariant controllers in which information propagates at least as fast as in the plant, convexity was established in~\cite{voubiabam03,bamvou05}. Similar algebraic characterization for a broader class of systems was introduced in~\cite{rotlal06}, and convexity was shown for problems with quadratically invariant constraint sets. Since these problems are convex in the impulse response parameters they are in general infinite dimensional. In~\cite{farjovAUT11}, a state-space description of systems in which information propagates at most one unit in space for every unit in time was provided and relaxations were used to obtain suboptimal controllers. In~\cite{twuege10}, the optimal control problem for switched autonomous systems was studied and optimality conditions for decentralization of multi-agent motions were derived. In~\cite{xiaboykim07}, convexity of the {\em symmetric\/} edge weight design for minimization of the mean-square deviation in distributed average consensus was shown.

While references~\cite{voubiabam03,bamvou05,rotlal06,farjovAUT11} focus on the design of optimal dynamic distributed controllers, we develop tools for the design of optimal static feedback gains with pre-specified structure. Even though the framework of~\cite{voubiabam03,bamvou05,rotlal06,farjovAUT11} does not apply to our setup, we identify a class of convex problems which can be cast as a semi-definite program (SDP). Furthermore, we show that the necessary conditions for optimality are given by coupled matrix equations, which can be solved by a combination of perturbation analysis and homotopy-based Newton's method. We consider the design of both symmetric and non-symmetric feedback gains and show that departure from optimal symmetric design can significantly improve the coherence of large-scale formations.

\subsection{Preview of key results}

We consider the design of optimal localized feedback gains for one-dimensional formations in which each vehicle only uses relative distances from its immediate neighbors and its own velocity. This nearest neighbor interaction imposes structural constraints on the feedback gains. We formulate the structured optimal control problem for both the single- and the double-integrator models. For single-integrators, we show that the structured optimal control problem is convex when we restrict the feedback gain to be a symmetric positive definite matrix. In this case, the global minimizer can be computed efficiently, and even analytical expressions can be derived. For double-integrators, we also identify a class of convex problems by restricting the controller to symmetric position and uniform diagonal velocity gains.

We then remove this symmetric restriction for both the single- and the double-integrator models and begin the design process with a spatially uniform controller. We develop a homotopy-based Newton's method that traces a continuous solution path from this controller to the optimal localized gain. Along this homotopy path, we solve a parameterized family of the structured optimal control problems and obtain {\em analytical\/} solutions when the homotopy parameter is small. We employ perturbation analysis to decouple the matrix equations that result from optimality conditions, thereby rendering the unique optimal structured gain. This solution is used to warm-start Newton's method in order to efficiently compute the desired optimal gains as the homotopy parameter is gradually increased.

In the second part of the paper, we examine how the performance of the optimally-controlled formation scales with the number of vehicles. We consider both macroscopic and microscopic performance measures based on whether attention is paid to the absolute position error of each vehicle or the relative position error between neighboring vehicles. We note that the macroscopic performance measure quantifies the resemblance of the formation to a rigid lattice, i.e., it determines the {\em coherence\/} of the formation. As shown in~\cite{bamjovmitpatTAC12}, even when local positions are well-regulated, an `accordion-like motion' of the formation can arise from poor scaling of the macroscopic performance measure (formation coherence) with the number of vehicles $N$. Our objective is thus to enhance formation coherence by means of optimal localized feedback design. In situations for which the control penalty in the quadratic performance objective is formation-size-independent we show that the optimal symmetric and non-symmetric controllers asymptotically provide $O(\sqrt{N})$ and $O(\sqrt[4]{N})$ scalings of formation coherence. However, this introduces similar growth of the control energy (per vehicle) with $N$. We show that bounded control energy can be obtained by judicious selection of an $N$-dependent control penalty, leading to $O(N)$ and $O(\sqrt{N})$ scalings of formation coherence for the optimal symmetric and non-symmetric controllers, respectively. These results are summarized in Table~\ref{tab.main} and they hold for both single- and double-integrators for formations in which each vehicle has access to {\em its own velocity\/}; see Sections~\ref{sec.scaling} and~\ref{sec.double} for additional details.

In addition to designing optimal localized controllers, we also provide an example of a spatially uniform non-symmetric controller that yields better scaling trends than the optimal spatially varying controller obtained by restricting design to symmetric gains. This indicates that departure from symmetry can improve coherence of large-scale formations and that the controller structure may play a more important role than the optimal selection of the feedback gains. On the other hand, our results also show that the optimal localized controller that achieves the best performance is both {\em non-symmetric and spatially-varying\/}.

If each vehicle has access to {\em its own velocity\/} and to relative distances from its nearest neighbors, we show similarity between the optimal position gains and performance scaling trends for single- and double-integrators. The latter observation is in agreement with analytical results obtained for spatially invariant formations~\cite{bamjovmitpatTAC12}. We note that performance of controllers that rely on relative measurements or unidirectional position exchange can differ significantly for these two models. For spatially-invariant formations with {\em relative position and velocity\/} measurements, it was shown in~\cite{bamjovmitpatTAC12} that the global performance scales as $O(N^3)$ for double-integrators and as $O(N)$ for single-integrators. In Section~\ref{sec.uni_nonsym_ctr}, we show that spatially uniform look-ahead strategy provides $O(\sqrt{N})$ scaling of the global performance for the single-integrator model. On the other hand, a look-ahead strategy that is not carefully designed can introduce unfavorable propagation of disturbances through formation of double-integrators~\cite{swahed96,seipanhed04}.

The paper is organized as follows. We formulate the structured optimal control problem in Section~\ref{sec.pro_for}, and show convexity of the symmetric gain design for the single-integrator model in Section~\ref{sec.sym}. For non-symmetric gains, we develop the homotopy-based Newton's method in Section~\ref{sec.hom_met}. We examine performance of localized controllers for the single- and the double-integrator models in Sections~\ref{sec.scaling} and~\ref{sec.double}, respectively, where we provide several explicit scaling relations. We conclude the paper in Section~\ref{sec.con_rem} with a brief summary of our contributions.

\section{Problem formulation}
    \label{sec.pro_for}

A system of $N$ identical vehicles moving along a straight line is shown in Fig.~\ref{fig.1Dplatoon}. All vehicles are equipped with ranging devices that allow them to measure relative distances with respect to their immediate neighbors. The objective is to design an {\em optimal\/} controller that uses {\em only local\/} information (i.e., relative distances between the neighboring vehicles) to keep each vehicle at its {\em global\/} position on a grid of regularly spaced points moving with a constant velocity.

We consider both the single- and the double-integrator models of the vehicles. The double-integrators are employed in many studies of vehicular formations; for example, see~\cite{levath66,melkuo71,swahed96,seipanhed04,jovbam05,
jovfowbamdan08,bamjovmitpatTAC12,barmehhes09,barmeh09}. On the other hand, the single-integrator (i.e., kinematic) model is simpler and perhaps more revealing in understanding the role of network topologies~\cite{jadlinmor03,olfmur04,marbrofra04,barhes06,xiaboykim07,zelmes11}. As we show in Section~\ref{sec.double}, the single- and the double-integrator models exhibit similar performance for formations in which each vehicle -- in addition to relative positions with respect to its immediate neighbors -- has an access to {\em its own velocity\/}. In the remainder of this section, we formulate the localized optimal control problem for both single- and double-integrators.

\subsection{Single- and double-integrator models}
    \label{sec.single}

We first consider the kinematic model in which control input $\bar{u}_n$ directly affects the velocity,
    \beq
    \dot{\bar{p}}_n
    \;=\;
    \bar{d}_n \,+\, \bar{u}_n,
    ~~
    n
    \, \in \,
    \{ 1,\ldots,N \},
    \non
    \eeq
where $\bar{p}_n$ is the position of the $n$th vehicle and $\bar{d}_n$ is the disturbance. The desired position of the $n$th vehicle is given by
    $
    p_{d,n}
    =
    v_d \, t + n \delta,
    $
where $v_d$ is the desired cruising velocity and $\delta$ is the desired distance between the neighboring vehicles. Every vehicle is assumed to have access to both $v_d$ and $\delta$. In addition, we confine our attention to formations with a known number of vehicles and leave issue of adaptation, merging, and splitting for future study.

The localized controller utilizes {\em relative position errors\/} between nearest neighbors,
    \[
    \bar{u}_n
    \,=\,
    \, - \,
    f_n (\bar{p}_n \,-\, \bar{p}_{n-1} \,-\, \delta)
    \,-\,
    b_n (\bar{p}_n \,-\, \bar{p}_{n+1} \,+\, \delta)
    \,+\,
    v_d,
    \]
where the design parameters $f_n$ and $b_n$ denote the forward and backward feedback gains of the $n$th vehicle. In deviation variables,
    $
    \{
    p_n
    :=
    \bar{p}_n - p_{d,n},
    $
    $u_n
    :=
    \bar{u}_n - v_d,
    $
    $
    d_n
    :=
    \bar{d}_n
    \},
    $
the single-integrator model with nearest neighbor interactions is given by
    \begin{subequations}
    \label{eq.single}
	\begin{align}
    \dot{p}_n
    \; = & \,~
    d_n \,+\, u_n,
    \label{eq.single-a}
    \\
    u_n
    \; = & \,~
    -
    f_n
    \left( p_n \,-\, p_{n-1} \right)
    \, - \,
    b_n
    \left( p_n \,-\, p_{n+1} \right),
    \label{eq.single-b}
    \end{align}
    \end{subequations}
where the relative position errors $p_n - p_{n-1}$ and $p_n - p_{n+1}$ can be obtained by ranging devices.

\begin{figure}[t]
    \begin{center}
    \subfloat[]{
    \includegraphics[width=0.45\textwidth]{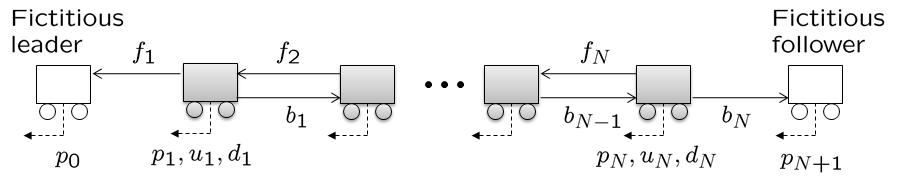}
    \label{fig.1D_fbgain_platoon}
    }
    \\
    \subfloat[]{
    \includegraphics[width=0.45\textwidth]{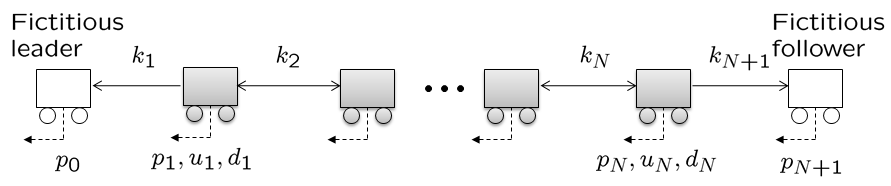}
    \label{fig.1D_sym_platoon}
    }
    \caption{Formation of vehicles with localized
    (a) non-symmetric; and
    (b) symmetric
    gains.}
    \end{center}
    \end{figure}

As illustrated in Fig.~\ref{fig.1D_fbgain_platoon}, fictitious lead and follow vehicles, respectively indexed by $0$ and $N+1$, are added to the formation. These two vehicles are assumed to move along their desired trajectories, implying that
    $
    p_0 = p_{N+1} = 0,
    $
and they are not considered to belong to the formation. Hence, the controls for the $1$st and the $N$th vehicles are given by
   \[
   \ba{rcl}
    u_1
    \!\!& = &\!\!
    - \, f_1 \, p_1
    \, - \,
    b_1
    \left( p_1 \,-\, p_2 \right),
    \\[0.1cm]
    u_N
    \!\!& = &\!\!
    - \, f_N
    \left( p_N \,-\, p_{N-1} \right)
    \, - \,
    b_N \, p_N.
   \ea
    \]
In other words, the first and the last vehicles have access to their own global position errors $p_1$ and $p_N$, which can be obtained by equipping them with GPS devices.

For the double-integrator model,
    \[
        \ddot{\bar{p}}_n
        \,=\,
        \bar{d}_n
        \,+\,
        \bar{u}_n,
        \;\;
        n
        \, \in \,
        \{ 1,\ldots,N \},
    \]
we consider the controller that has an access to the relative position errors between the neighboring vehicles and the absolute velocity errors,
    \[
    \ba{rcl}
        \bar{u}_n
        \!\!& = &\!\!
        - \, f_n \, (\bar{p}_n \,-\, \bar{p}_{n-1} \,-\, \delta)
        \,-\,
        b_n \, (\bar{p}_n \,-\, \bar{p}_{n+1} \,+\, \delta)
        \\[0.1cm]
        \!\!& &\!\!
        -\,
        g_n \, (\dot{\bar{p}}_n \,-\, v_d),
    \ea
    \]
where $g_n$ denotes the velocity feedback gain. In deviation variables,
    $
    \{
    p_n
    :=
    \bar{p}_n - p_{d,n}
    $,
    $
    v_n
    :=
    \dot{\bar{p}}_n - v_d,
    $
    $
    u_n
    :=
    \bar{u}_n,
    $
    $
    d_n
    :=
    \bar{d}_n
    \},
    $
the double-integrator model is given by
    \begin{subequations}
    \label{eq.double}
	\begin{align}
    \ddot{p}_n
    \; = & \,~
    d_n \,+\, u_n,
    \label{eq.double-a}
    \\
    u_n
    \; = & \,~
    -
    f_n \left( p_n \,-\, p_{n-1} \right)
    \, - \,
    b_n \left( p_n \,-\, p_{n+1} \right)
    \, - \,
    g_n \, v_n.
    \label{eq.double-b}
    \end{align}
    \end{subequations}

In matrix form, control laws~(\ref{eq.single-b}) and~(\ref{eq.double-b}) can be written as,
    \[
    \ba{rcl}
    u
    & \! = \! &
    - \, F C \, p
    \,=\,
    - \, \obt{F_f}{F_b}
    \left[\ba{c} C_f \\ C_f^T \ea \right] p,
    \\[0.3cm]
    u
    & \! = \! &
    - \, F C \, \tbo{p}{v}
    \\
    & \! = \! &
    - \,
    \obth{F_f}{F_b}{F_v}
    \left[ \ba{cc} C_f & O\\ C_f^T & O \\ O & I \ea \right]
    \tbo{p}{v},
    \ea
    \]
where $p$, $v$, and $u$ denote the position error, the velocity error, and the control input vectors, e.g.,
    $
    p =
    \left[\, p_1 ~ \cdots ~ p_N \,\right]^T.
    $
Furthermore, the $N \times N$ diagonal feedback gains are determined by
    \[
    F_f
    \, := \,
    \diag \, \{f_n\},
    ~~
    F_b
    \, := \,
    \diag \, \{b_n\},
    ~~
    F_v
    \, := \,
    \diag \, \{g_n\},
    \]
and $C_f$ is a sparse Toeplitz matrix with $1$ on the main diagonal and $-1$ on the first lower sub-diagonal. For example, for $N = 4$,
    \beq
    F_f
    \, = \,
    \left[
    \ba{cccc}
    \!f_1\! & 0 & 0 & 0 \\
    0 & \!f_2\! & 0 & 0 \\
    0 & 0 & \!f_3\! & 0 \\
    0 & 0 & 0 & \!f_4\! \\
    \ea
    \right],
    ~~
    C_f
    \, = \,
    \left[
    \ba{rrrc}
    1 & 0 & 0 & 0 \\
    \!-1 & 1 & 0 & 0 \\
    0 & \!-1 & 1 & 0 \\
    0 & 0 & \!-1 & 1
    \ea
    \right].
    \label{eq.Cf}
    \eeq
Thus, $C_f \, p$ determines the vector of the relative position errors $p_n - p_{n-1}$ between each vehicle and the one in front of it; similarly, $C_f^T \, p$ determines the vector of the relative position errors $p_n - p_{n+1}$ between each vehicle and the one behind it.

We will also consider formations with no fictitious followers. In this case, the $N$th vehicle only uses relative position error with respect to the $(N-1)$th vehicle, i.e., $b_N = 0$ implying that $u_N = -f_N \, (p_N - p_{N-1})$ for the single-integrator model and $u_N = -f_N \, (p_N - p_{N-1}) - g_N \, v_N $ for the double-integrator model.

\subsection{Structured ${\cal H}_2$ problem}

The state-space representation of the vehicular formation is given by
    \beq
    \ba{rcl}
    \dot{x}
    &=&
    A \, x
    \,+\,
    B_1 \, d
    \,+\,
    B_2 \, u,
    \\
    y
    &=&
    C \, x,
    ~~
    u
    ~=\,
    - \, F \, y.
    \ea
    \tag*{(SS)}
    \label{SS}
    \eeq
For the single-integrator model~(\ref{eq.single}), the state vector is $x = p$, the measured output $y$ is given by the relative position errors between the neighboring vehicles, and
    \beq
    \ba{l}
    A
    ~=~
    O,
    ~~
    B_1
    \, = \,
    B_2
    \, = \,
    I,
    ~~
    C
    ~=~
    \left[
    \ba{c}
    C_f
    \\
    C_f^T
    \ea
    \right],
    \\
    F
    \, = \,
    \left[
    \ba{cc}
    F_f & F_b
    \ea
    \right].
    \ea
    \tag*{(VP1)}
    \label{VP1}
    \eeq
For the double-integrator model~(\ref{eq.double}), the state vector is $x = [\, p^T ~ v^T \,]^T$, the measured output $y$ is given by the relative position errors between the neighboring vehicles and the absolute velocity errors, and
    \beq
    \ba{l}
    \tag*{(VP2)}
    \label{VP2}
    A \,=\, \tbt{O}{I}{O}{O},
    ~
    B_1 \,= \, B_2 \,=\, \tbo{O}{I},
    \\[0.25cm]
    C \,=\, \left[ \ba{cc} C_f & O\\ C_f^T & O \\ O & I \ea \right],
    ~
    F \,=\, \obth{F_f}{F_b}{F_v}.
    \ea
    \eeq
Here, $O$ and $I$ denote the zero and identity matrices, and $\{ F_f$, $F_b$, $F_v$, $C_f \}$ are defined in~(\ref{eq.Cf}).

Upon closing the loop, we have
    \beq \nn
    \ba{rcl}
    \dot{x}
    & \!\! = \!\! &
    (A \,-\, B_2 F C) \, x \,+\, B_1 d,
    \\[0.2cm]
    z
    & \!\! = \!\! &
    \tbo{Q^{1/2} \, x}{r^{1/2} \, u}
    \,=\,
    \tbo
    {Q^{1/2}}
    {- r^{1/2} F C}
    x,
    \ea
    \eeq
where $z$ encompasses the penalty on both the state and the control. Here, $Q$ is a symmetric positive semi-definite matrix and $r$ is a positive scalar. The objective is to design the {\em structured feedback gain\/} $F$ such that the influence of the white stochastic disturbance $d$, with zero mean and unit variance, on the performance output $z$ is minimized (in the ${\cal H}_2$ sense). This control problem can be formulated as~\cite{farlinjovCDC09,linfarjovTAC11AL}
    \beq
    \tag*{(SH2)}
    \label{SH2}
    \ba{ll}
    \!\!\!\!
    \text{minimize}
    & \!\!
    J
    \,=\,
    \trace \left( P B_1 B_1^T \right)
    \\[0.15cm]
    \!\!\!\!
    \text{subject to}
    & \!\!
        (A \,-\, B_2 F C)^T P \,+\, P (A \,-\, B_2 F C)
        \; =
    \\[0.1cm]
    &
        ~~~~~
        - \, (Q \,+\,  r \, C^T F^T F C),
        ~~
        F \in \eS
    \ea
    \eeq
where $\eS$ denotes the structural subspace that $F$ belongs to.

As shown in~\cite{farlinjovCDC09}, the necessary conditions for
optimality of \ref{SH2} are given by the set of coupled matrix
equations in $F$, $P$, and $L$
    \begin{align*}
    \!\!\!
    & (A \,-\, B_2 F C)^T \, P \,+\, P \, (A \,-\, B_2 F C)
    \,=\,
    \\
    &
    \hspace{3.5cm}
    - \left( Q \,+\,  r \, C^T F^T F C \right),
    \tag*{(NC1)}
    \label{NC1}
    \\
    \!\!\!
    & (A \,-\, B_2 F C) \, L \,+\, L \, (A \,-\, B_2 F C)^T  \,=\,
    - \, B_1 B_1^T,
    \tag*{(NC2)}
    \label{NC2}
    \\
    \!\!\!
    &
    ( r F C L C^T) \circ I_\eS
    \, = \,
    (B_2^T P L C^T) \circ I_\eS.
    \tag*{(NC3)}
    \label{NC3}
    \end{align*}
Here, $P$ and $L$ are the closed-loop observability and controllability Gramians, $\circ$ denotes the entry-wise multiplication of two matrices, and the matrix $I_\eS$ in~\ref{NC3} denotes the structural identity of the subspace $\eS$ under the entry-wise multiplication, i.e.,
    $
    F \circ I_\eS
    =
    F,
    $
with
    $
        I_{\eS} = [\,I \;\; I \,]
    $
for the single-integrator model and
    $
        I_{\eS} = [\,I \;\; I \;\; I \,]
    $
for the double-integrator model. (For example, $[\,F_f \;\; F_b \,] \circ [\,I \;\; I \,] = [\,F_f \;\; F_b \,]$.) In the absence of the fictitious follower, an additional constraint $b_N = 0$ is imposed in~\ref{SH2} and thus, the structural identity for the single- and the double-integrator models are given by $ [\, I \;\; I_z \,] $ and $[\, I \;\; I_z \;\; I \,]$, respectively. Here, $I_z$ is a diagonal matrix with its main diagonal given by $[\,1 \,\cdots\, 1 \,\, 0\,]$.

\remark
Throughout the paper, the structured optimal feedback gain $F$ is obtained by solving~\ref{SH2} with $Q = I$. This choice of $Q$ is motivated by our desire to enhance formation coherence, i.e., to keep the global position and velocity errors $p_n$ and $v_n$ small using {\em localized feedback}. Since the methods developed in the paper can be applied to other choices of $Q$, we will describe them for general $Q$ and set $Q = I$ when presenting computational results.

\subsection{Performance of optimal localized controller}
    \label{sec.Pi}

To evaluate the performance of the optimal localized controller $F$, obtained by solving~\ref{SH2} with $Q = I$, we consider the closed-loop system
    \beq
    \label{eq.CL}
    \ba{rcl}
    \dot{x}
    & \!\! = \!\! &
    (A \,-\, B_2 F C) \, x \,+\, B_1 d,
    \\[0.2cm]
    \zeta
    & \!\! = \!\! &
    \tbo{\zeta_1}{\zeta_2}
    \,=\,
    \tbo{Q_s^{1/2}}{ - \, F C } x,
    ~~
    s \, = \, g
    ~~
    \mbox{or}
    ~~
    s \, = \, l,
    \ea
    \eeq
where $\zeta_1$ is the global or local performance output and $\zeta_2$ is the control input. Motivated by~\cite{bamjovmitpatTAC12}, we examine two state performance weights for the single-integrator model
    \bi
    \item
    Macroscopic (global):  $Q_g \, = \, I$;
    \item
    Microscopic (local): $Q_l \, = \, T$,
    \ei
where $T$ is an $N \times N$ symmetric Toeplitz matrix with its first row given by
$\left[ \, 2 ~ -1 ~\, 0 ~ \cdots ~ 0 \, \right] \in \bbR^{N}$.
For example, for $N = 4$,
    \beq \label{eq.T}
    T
    \, = \,
    \left[
    \ba{rrrr}
    2    &   -1  & 0    & 0 \\
    -1   &   2   & -1   & 0 \\
    0    &   -1  & 2    & -1 \\
    0    &   0   & -1   & 2
    \ea
    \right].
    \eeq
The macroscopic performance weight $Q_g = I$ penalizes the {\em global\/} (absolute) position errors,
    \[
        \zeta_1^T \zeta_1
        \,=\,
        p^T Q_g \, p
        \,=\,
        \sum_{n\,=\,1}^N p_n^2,
    \]
and the microscopic performance weight $Q_l = T$ penalizes the
{\em local\/} (relative) position errors,
    \[
        \zeta_1^T \zeta_1
        \,=\,
        p^T Q_l \, p
        \,=\,
        \sum_{n\,=\,0}^N (p_{n}-p_{n+1})^2,
    \]
with $p_0 = p_{N+1} = 0.$ These state weights induce the macroscopic and microscopic performance measures~\cite{bamjovmitpatTAC12} determined by the {\em formation-size-normalized\/} ${\cal H}_2$ norm
    \beq
    \Pi_s (N)
    \, = \,
    (1/N) \, \| G_1 \|_2^2,
    ~~
    s \, = \, g
    ~~
    \mbox{or}
    ~~
    s \, = \, l,
    \tag*{($\Pi$)}
    \label{Pi}
    \eeq
where $G_1$ is the transfer function of~(\ref{eq.CL}) from $d$ to $\zeta_1$. The macroscopic performance measure $\Pi_g$ quantifies the resemblance of the formation to a rigid lattice, i.e., it determines the {\em coherence\/} of the formation~\cite{bamjovmitpatTAC12}. On the other hand, the microscopic performance measure $\Pi_l$ quantifies how well regulated the distances between the neighboring vehicles are. We will also examine the formation-size-normalized control energy (variance) of the closed-loop system~(\ref{eq.CL}),
    \beq
    \Pctr (N)
    \, = \,
    (1/N) \, \| G_2 \|_2^2,
    \non
    \eeq
which is determined by the ${\cal H}_2$ norm of the transfer function $G_2$ from $d$ to $\zeta_2 = u$.

Similarly, for the double-integrator model, we use the following performance weights
    \bi
    \item
    Macroscopic (global),  $Q_g \, = \, \tbt{I}{O}{O}{I}$;
    \\[0.1cm]
    \item
    Microscopic (local), $Q_l \, = \, \tbt{T}{O}{O}{I}$.
    \ei

\subsection{Closed-loop stability: the role of fictitious vehicles}
    \label{sec.stability_nofollower}

We next show that at least one fictitious vehicle is needed in order to achieve closed-loop stability. This is because the absence of GPS devices in the formation prevents vehicles from tracking their absolute desired trajectories.

For the single-integrator model, the state-feedback gain
    $
    K_p
    =
    F_f \, C_f
    +
    F_b \, C_f^T
    $
is a structured tridiagonal matrix satisfying
    $
        K_p \, \bfo
        =
        \big[ \, f_1 \;\; 0 \;\; \cdots \;\; 0 \;\; b_N \, \big]^T
    $
where $\bfo$ is the vector of all $1$'s. If neither the $1$st nor the $N$th vehicle has access to its own global position, i.e., $f_1 = b_N = 0$, then $K_p$ has a zero eigenvalue with corresponding eigenvector $\bfo$. Hence, the closed-loop system is not asymptotically stable regardless of the choice of the feedback gains $\{f_n\}_{n\,=\,2}^N$ and $\{b_n\}_{n\,=\,1}^{N-1}$. In the presence of stochastic disturbances, the average-mode (associated with the eigenvector $\bfo$) undergoes a random walk and the steady-state variance of the deviation from the absolute desired trajectory becomes unbounded~\cite{xiaboykim07,bamjovmitpatTAC12,zelmes11}. In this case, other performance measures that render this average-mode {\em unobservable\/} can be considered~\cite{bamjovmitpatTAC12}.

For the double-integrator model, the action of
    $
        A_{\rm cl}
        =
        A
        -
        B_2 F C
    $
on  $[ \, \bfo^T \;\; \mathbf{0}^T \,]^T$ is given by
    \[
    \tbt{O}{I}{- \, K_p}{- \, F_v} \tbo{\bfo}{\bf 0}
    \,=\,
    \tbo{\bf 0}{- \, K_p \, \bfo} ,
    \]
where ${\bf 0}$ is the $N$-vector of all $0$'s. Thus, if $f_1 = b_N = 0$ then $A_{\rm cl}$ has a zero eigenvalue with corresponding eigenvector $[ \, {\bf 1}^T \;\; {\bf 0}^T \,]^T$. Therefore, for both the single- and the double-integrator models, we need at least one vehicle with access to its global position in order to achieve closed-loop stability.

\section{Design of symmetric gains for the single-integrator model: a convex problem}
    \label{sec.sym}

In this section, we design the optimal symmetric feedback gains for the single-integrator model; see Fig.~\ref{fig.1D_sym_platoon}. This is a special case of the localized design, obtained by restricting the forward and the backward gains between the neighboring vehicles to be equal to each other, i.e., $f_n=b_{n-1}$ for $n \in \{2,\ldots,N\}$. Under this assumption, we show that~\ref{SH2} is a convex optimization problem for the single-integrator model. This implies that the {\em global minimum\/} can be computed efficiently. Furthermore, in the absence of the fictitious follower, we provide {\em analytical expressions\/} for the optimal feedback gains.

Let us denote $k_1=f_1$ and $k_{N+1}=b_N$ and let
    \beq
    \label{eq.sym_gain}
    k_n \,=\, f_n \,=\, b_{n-1},
    ~~
    n \in  \{2, \ldots, N \}.
    \eeq
For the single-integrator model, the structured gain becomes a {\em symmetric\/} tridiagonal matrix
    \beq
    \label{eq.K}
    \ba{rcl}
        K
        \!\!& = &\!\!
        F_f \, C_f
        \; + \;
        F_b \, C_f^T
        \\
        \!\! & = & \!\!
        \left[
        \ba{ccccc}
        k_1+k_2 & -k_2    &         &       \\
        -k_2    & k_2+k_3 & \ddots  &       \\
                & \ddots  & \ddots  & -k_N   \\
                &         & -k_N    & k_N+k_{N+1}
        \ea
        \right].
        \ea
    \eeq
Consequently, $A_{\rm cl} = -K$ is Hurwitz if and only if $K$ is positive
definite, in which case the Lyapunov equation in~\ref{SH2} simplifies to
    \[
        K P \,+\, P K \,=\, Q \,+\, r K K.
    \]
The application of~\cite[Lemma~1]{bamdah03} transforms the problem~\ref{SH2} of optimal symmetric design for the single-integrator model to
    \beq
    \nn
    \tag*{(SG)}
    \label{SG}
    \ba{ll}
    \!\!
    \underset{K}{\operatorname{minimize}}
    & \!\!
    J(K)
    \; = \;
    (1/2) \, \trace \left( Q K^{-1} \,+\, r K \right)
    \\[0.15cm]
    \!\!
    \text{subject to}
    & \!\!
    K \, > \, 0 ~\text{and $K \in \eS_K$}
    \ea
    \eeq
where $K \in \eS_K$ is a linear structural constraint given by (\ref{eq.K}). (Specifically, $K = F_f C_f + F_b C_f^T$ is a symmetric tridiagonal matrix with the linear constraint~(\ref{eq.sym_gain}).) By introducing an auxiliary variable $X = X^T \geq Q^{1/2} K^{-1} Q^{1/2}$, we can formulate~\ref{SG} as an SDP in $X$ and $K$
    \beq
    \ba{ll}
    \!\!\!\!
    \underset{X, \, K}{\operatorname{minimize}}
    &
    (1/2)
    \, \trace \left( X \, + \, r K \right)
    \\[0.cm]
    \!\!\!\!
    \text{subject to}
    &
    K \, > \, 0,
    ~~
    K \in \eS_K,
    ~~
    \left[\ba{cc}
    K & Q^{1/2} \\
    Q^{1/2} & X
    \ea\right] \geq 0,
    \ea
    \!\!\!
    \nn
    \eeq
which can be solved using available SDP solvers. Here, we have used the Schur complement~\cite[Appendix~A.5.5]{boyvan04} in conjunction with $K > 0$ to express $X \geq Q^{1/2} K^{-1} Q^{1/2}$ as an LMI.

Next, we exploit the structure of $K$ to express $J$ in~\ref{SG} with $Q = I$ in terms of the feedback gains $\{ k_n \}_{n \,=\, 1}^{N + 1}$ between the neighboring vehicles. Since the inverse of the symmetric tridiagonal matrix $K$ can be determined analytically~\cite[Theorem 2.3]{meu92}, the $ij$th entry of $K^{-1}$ is given by
    \beq
    (K^{-1})_{ij}
    \, = \,
    \dfrac{\gamma_i \, (\gamma_{N+1}-\gamma_j)}{\gamma_{N+1}},
    ~~
    j \geq i,
    ~~~
    \gamma_i
    \, = \,
    \ds{\sum_{n\,=\,1}^i} \frac{1}{k_n},
    \label{eq.invK}
    \eeq
yielding the following expression for $J$
    \begin{align*}
        J
        \, &= \,
        (1/2)
        \,
        \trace
        \left(
        K^{-1}
        \,+\,
        r K
        \right)
        \\
        \, &= \,
        \dfrac{1}{2} \sum_{n\,=\,1}^N
        \frac{\gamma_n(\gamma_{N+1}-\gamma_n)}{\gamma_{N+1}}
        \, + \,
        r
        \left(
        \dfrac{k_1 + k_{N+1}}{2}
        \, + \,
        \sum_{n\,=\,2}^N k_n
        \right).
    \end{align*}
The above expression for $J$ is well-defined for $\{ k_n \}_{n \,=\, 1}^{N + 1}$ that guarantee positive definiteness of $K$ in~(\ref{eq.K}); this is because the closed-loop $A$-matrix is determined by $A_{\rm cl} = -K$. The global minimizer of $J$ can be computed using the gradient method; see Appendix~\ref{app.grad_sym_K}.

For the formations without the fictitious follower, we next derive {\em explicit analytical expression\/} for the global symmetric minimizer $K = K^T > 0$ of~\ref{SG} with $Q = I$. In this case $k_{N+1} = 0$ and the $ij$th entry of $K^{-1}$ in~(\ref{eq.invK}) simplifies to $(K^{-1})_{ij} = \gamma_i$ for $ j \geq i $. Consequently, the unique minimum of
    \begin{align*}
        J
        \;&=\;
        \dfrac{1}{2}
        \sum_{n\,=\,1}^N \gamma_n
        \,+\,
        r
        \left(
        \dfrac{k_1}{2}
        \,+\,
        \sum_{n\,=\,2}^N k_n
        \right)
        \\
        \;&=\;
        \dfrac{1}{2}
        \sum_{n\,=\,1}^N
        \dfrac{N+1-n}{k_n}
        \,+\,
        r
        \left(
        \dfrac{k_1}{2}
        \,+\,
        \sum_{n\,=\,2}^N k_n
        \right),
    \end{align*}
is attained for
    \beq \label{eq.symK_QI}
        k_1
        \, = \,
        \sqrt{N/r},
        ~~
        k_n
        \, = \,
        \sqrt{(N+1-n)/(2 r)},
        ~~
        n
        \in
        \{ 2,\ldots,N \}.
    \eeq
We also note that
    \beq
    \label{eq.KequKinv}
    \ba{rcl}
    \trace
    \left(
    K^{-1}
    \right)
    \!\! & = & \!\!
    \ds{\sum_{n\,=\,1}^N \gamma_n}
    \;=\;
    \ds{\sum_{n\,=\,1}^N}
    \dfrac{N+1-n}{k_n}
    \\[0.1cm]
    \!\! & = & \!\!
    r
    \left(
    k_1
    \,+\,
    2
    \ds{\sum_{n\,=\,2}^N}
    k_n
    \right)
    \,=\,
    r
    \,
    \trace
    \,
    (
    K
    ),
    \ea
    \eeq
where the third equality follows from~(\ref{eq.symK_QI}). This result is used to examine the performance of large-scale formations in Section~\ref{sec.opt_ctr}.

\begin{figure}[h]
    \centering
    \includegraphics[width=0.35\textwidth]{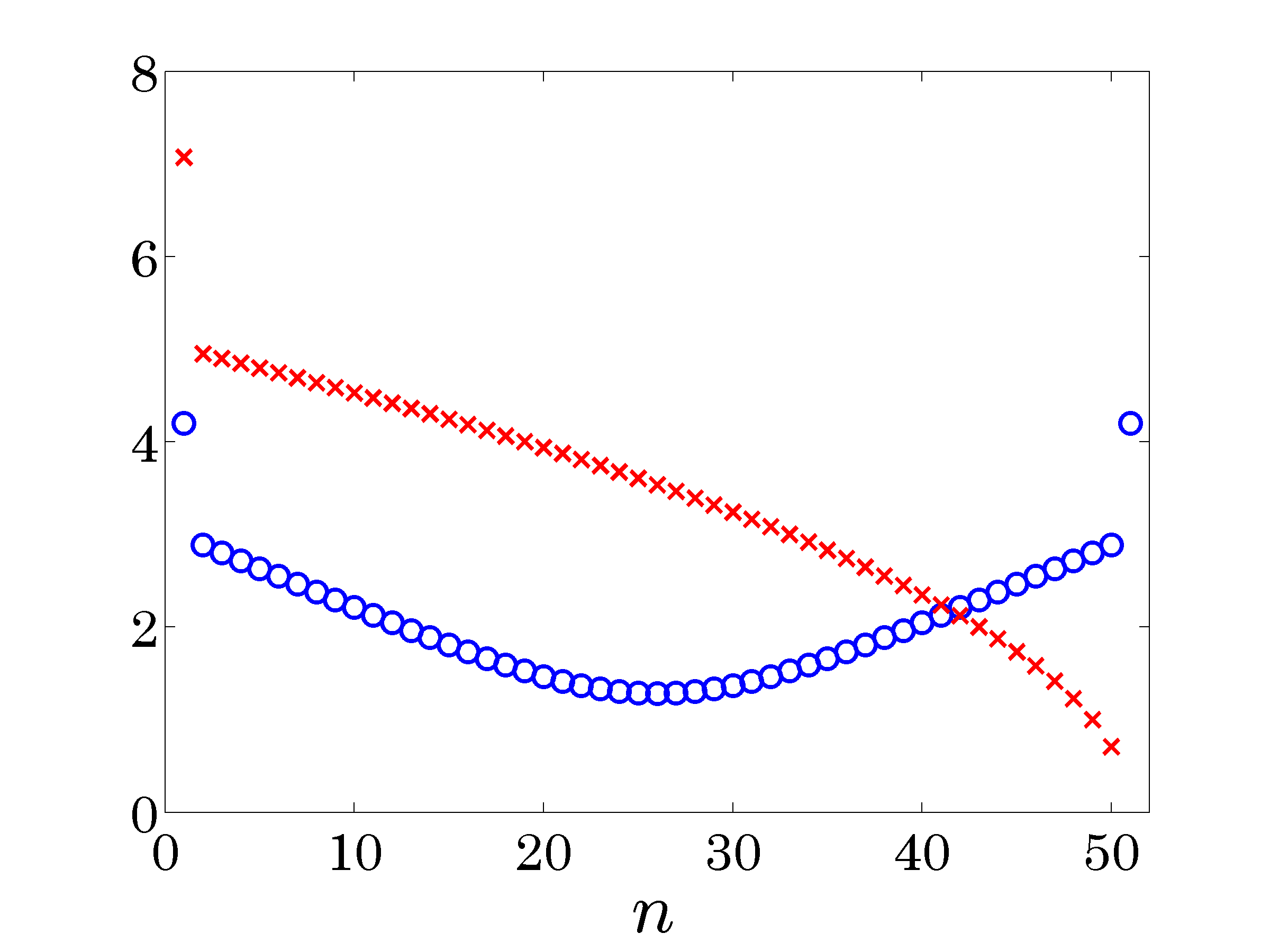}
    \caption{Optimal symmetric gains for formations with follower ($\circ$) and without follower ($\times$) for $N = 50$, $Q = I$, and $r = 1$.
    ($\times$) are obtained by evaluating formula~(\ref{eq.symK_QI}) and ($\circ$) are computed using the gradient method described in Appendix~\ref{app.grad_sym_K}.}
    \label{fig.sym_gain}
    \end{figure}

Figure~\ref{fig.sym_gain} shows the optimal symmetric gains for a formation with $N = 50$ vehicles, $Q = I$, and $r = 1$. Since the fictitious leader and the follower always move along their desired trajectories, the vehicles that are close to them have larger gains than the other vehicles. When no fictitious follower is present, the gains decrease monotonically from the first to the last vehicle; see $(\times)$ in Fig.~\ref{fig.sym_gain}. In other words, the farther away the vehicle is from the fictitious leader the less weight it places on the information coming from its neighbors. This is because uncorrelated disturbances that act on the vehicles corrupt the information about the absolute desired trajectory as it propagates from the fictitious leader down the formation (via relative information exchange between the vehicles). When both the fictitious leader and the follower are present, the gains decrease as one moves from the boundary to the center of the formation; see $(\circ)$ in Fig.~\ref{fig.sym_gain}. This can be attributed to the fact that the information about the absolute desired trajectories becomes noisier as it propagates from the fictitious vehicles to the center of the formation.

\section{Homotopy-based Newton's method}
    \label{sec.hom_met}

In this section, we remove the symmetric feedback gain restriction and utilize
a homotopy-based Newton's method to solve~\ref{SH2}. In~\cite{farlinjovCDC09}, Newton's method for general structured ${\cal H}_2$ problems is developed.
For~\ref{SH2} with the specific problem data \ref{VP1} and \ref{VP2}, it is possible to employ a homotopy-based approach to solve a parameterized family of problems, which ranges between an easily solvable problem and the problem of interest.

In particular, we consider
    \beq
    \label{eq.Q}
        Q(\veps)
        \; = \;
        Q_0
        \; + \;
        \veps
        \left(
        Q_d \,-\, Q_0
        \right),
    \eeq
where $Q_0$ is the initial weight to be selected, $Q_d$ is the {\em desired\/} weight, and $\veps \in [0,1]$ is the homotopy parameter. Note that $Q = Q_0$ for $\veps = 0$, and $Q = Q_d$ for $\veps = 1$. The homotopy-based Newton's method consists of three steps: (i) For $\veps = 0$, we find the initial weight $Q_0$ with respect to which a spatially uniform gain $F_0$ is {\em inversely optimal\/}. This is equivalent to solving problem \ref{SH2} {\em analytically\/} with the performance weight $Q_0$. (ii) For $0 < \veps \ll 1$, we employ perturbation analysis to determine the first few terms in the expansion
    $
        F(\veps) = \sum_{n\,=\,0}^\infty \veps^n F_n.
    $
(iii) For larger values of $\veps$, we use Newton's method for structured ${\cal H}_2$ design~\cite{farlinjovCDC09} to solve \ref{SH2}. We gradually increase $\veps$ and use the structured optimal gain obtained for the previous value of $\veps$ to initialize the next round of iterations. This process is repeated until the desired value $\veps = 1$ is reached.

In the remainder of this section, we focus on the single-integrator model. In Section~\ref{sec.double}, we solve problem~\ref{SH2} for the double-integrator model.

\subsection{Spatially uniform symmetric gain: inverse optimality for $\veps=0$}
    \label{sec.inv_opt_design}

One of the simplest localized strategies is to use {\em spatially uniform gain\/}, where $F_f$ and $F_b$ are diagonal matrices with $f_n = f$ and $b_n = b$ for all $n$ and some positive $f$ and $b$. In particular, for $F_f = F_b = I$ it is easy to show closed-loop stability and to find the performance weight $Q_0$ with respect to which the spatially uniform symmetric gain
    \[
    K_0
    \, = \,
    F_0 \, C
    \, = \,
    \left[
    \ba{cc}
    I & I
    \ea
    \right]
    \left[
    \ba{c}
    C_f
    \\
    C_f^T
    \ea
    \right]
    \, = \,
    T
    \]
is inversely optimal. The problem of inverse optimality amounts to finding the performance weight $Q_0$ for which an {\em a priori\/} specified $K_0$ is the corresponding optimal state-feedback gain~\cite{kal64,jovIO10}. From linear quadratic regulator theory, the optimal state-feedback gain is given by
    $
    K_0
    =
    R^{-1} B_2^T P_0
    $
where $P_0$ is the positive definite solution of
    \beq \nn
        A^T P_0 \,+\, P_0 A \,+\, Q_0 \,-\, P_0 B_2 R^{-1} B_2^T P_0
        \,=\,
        0.
    \eeq
For the kinematic model~\ref{VP1}, $A = O$ and $B_2 = I$, with $R = r I$, we have
    $ K_0 = r^{-1} P_0  $
and
    $ Q_0 - r^{-1} P_0 P_0   =  0 $.
Therefore, the state penalty
    $
    Q_0
    =
    r K_0^2
    =
    r T^2
    $
guarantees inverse optimality of the spatially uniform symmetric gain $K_0$.
The above procedure of finding $Q_0$ can be applied to any structured gain $F_0$ that yields a symmetric positive definite $K_0$, e.g., the optimal symmetric gain of Section~\ref{sec.sym}.

\subsection{Perturbation analysis for $\veps \ll 1$}
    \label{sec.pert_ana}

We next utilize perturbation analysis to solve~\ref{SH2}
with $Q(\veps$) given by (\ref{eq.Q}) for $\veps \ll 1$. For small $\veps$, by representing $P$, $L$, and $F$ as
    \begin{align*}
    P \,= \sum_{n \,=\, 0}^{\infty} \veps^n P_n,
    ~~
    L \,= \sum_{n \,=\, 0}^{\infty} \veps^n L_n,
    ~~
    F \,= \sum_{n \,=\, 0}^{\infty} \veps^n F_n,
    \end{align*}
substituting in \ref{NC1}-\ref{NC3}, and collecting same-order terms in $\veps$, we obtain the set of equations~\ref{PA} with
    $
    A_0
    :=
    A
    -
    B_2 F_0 C.
    $
Note that these equations are {\em conveniently coupled\/} in one direction, in the sense that for any $n \geq 1$, $O(\veps^{n})$ equations depend only on the solutions of $O(\veps^m)$ equations for $m \leq n$. In particular, it is easy to verify that the first and the third equations of $O(1)$ are satisfied with
    $
    K_0
    =
    F_0 C
    =
    r^{-1} B_2^T P_0
    $
and with $Q_0 = r K_0^2$ identified in Section~\ref{sec.inv_opt_design}. Thus, the matrix $L_0$ can be obtained by solving the second equation of $O(1)$, and the matrices $P_1$, $F_1$, and $L_1$ can be obtained by solving the first, the third, and the second equations of $O(\veps)$, respectively. The higher order terms $F_n$, $P_n$, and $L_n$ can be determined in a similar fashion. The matrix $F$ found by this procedure is the {\em unique optimal\/} solution of the control problem~\ref{SH2} for $\veps \ll 1$. This is because the equations~\ref{PA}, under the assumption of convergence for small $\veps$, give a unique matrix $F(\veps) = \sum_{n \,=\,0}^{\infty} \veps^n F_n$.

    \begin{figure*}
    \[
    \ba{l}
    O(1):
    \left\{
    \ba{rcl}
    A_0^T P_0 \,+\, P_0 A_0
    & \!\! = \!\! &
    -(Q_0 \,+\, r \, C^T F_0^T F_0 C)\\[0.1cm]
    A_0 L_0 \,+\, L_0 A_0^T
    & \!\! = \!\! &
    - B_1 B_1^T \\[0.1cm]
    (r F_0 C L_0 C^T) \circ I_\eS
    & \!\! = \!\! &
    (B_2^T P_0 L_0 C^T) \circ I_\eS
    \ea
    \right.
    \\[0.75cm]
    O(\veps):
    \left\{
    \ba{rcl}
    A_0^T P_1 \,+\, P_1 A_0
    & \!\! = \!\! &
    - (Q_d - Q_0 )\\[0.1cm]
    A_0 L_1 \,+\, L_1 A_0^T
    & \!\! = \!\! &
    (B_2 F_1 C) L_0 \,+\, L_0 (B_2 F_1 C)^T\\[0.1cm]
    (r F_1 C L_0 C^T) \circ I_\eS
    & \!\! = \!\! &
    (B_2^T P_1 L_0 C^T) \circ I_\eS
    \ea
    \right.
    \\[0.75cm]
    O(\veps^2):
    \left\{
    \ba{rcl}
    A_0^T P_2 \,+\, P_2 A_0
    & \!\! = \!\! &
    (B_2 F_1 C)^T P_1 \,+\, P_1 (B_2 F_1 C) \,-\, r C^T F_1^T F_1 C \\[0.1cm]
    A_0 L_2 \,+\, L_2 A_0^T
    & \!\! = \!\! &
    (B_2 F_1 C) L_1 \,+\, L_1 (B_2 F_1 C)^T \,+\, (B_2 F_2 C)L_0 \,+\, L_0(B_2 F_2 C)^T\\[0.1cm]
    (r F_2 C L_0 C^T ) \circ I_\eS
    & \!\! = \!\! &
    (B_2^T P_1 L_1 C^T \,+\, B_2^T P_2 L_0 C^T \,-\, r F_1 C L_1 C^T ) \circ I_\eS
    \ea
    \right.
    \\[0.35cm]
    ~~~\vdots
    \hspace{4.1cm}
    \left.
    \ba{rcl}
    ~
    & \!\! \vdots \!\! &
    ~
    \ea
    \right.
    \ea
    \tag*{(PA)}
    \label{PA}
    \]
\hrulefill
\end{figure*}

We next provide analytical expressions for $F_1 = [ \, F_f^{(1)} ~\, F_b^{(1)} \, ]$ obtained by solving the $O(\veps)$ equations in~\ref{PA} with $r = 1$, $Q_0 = T^2$, and $Q_d = I$. When a fictitious follower is present, we have (derivations are omitted for brevity)
    \beq
    \label{eq.F1}
    \ba{l}
    f_n^{(1)}
     =
    \dfrac{n (n-N-1) (4 n (N+1)-N (2 N+7)+1)}{12 \left(N^2-1\right)}
     -
    \dfrac{1}{2},
    \\[0.25cm]
    b_n^{(1)}
     =
    \dfrac{n (N + 1 - n) (4 n (N+1)-N (2 N+1)-5)}{12 \left(N^2-1\right)}
     -
    \dfrac{1}{2},
    \ea
    \eeq
where $f_n^{(1)}$ and $b_n^{(1)}$ denote the $n$th diagonal entries of $F_f^{(1)}$ and $F_b^{(1)}$. From~(\ref{eq.F1}) it follows that
    $
    f_n^{(1)}
    =
    b_{N + 1 - n}^{(1)}
    $
    for
    $
    n \in \{1, \ldots, N\}.
    $
When a fictitious follower is not present, we have
    \[
    \ba{ll}
    f_n^{(1)}
    \; = \;
    (- \, n^2 \,+\, (N+1)n \,-\, 1)/2,
    &\!\!
    n \in \{1,\ldots,N - 1\},
    \\
    f_N^{(1)} \,=\, (N - 1)/2,
    &
    \\
    b_n^{(1)}
    \; = \;
    (n^2 \,-\, Nn \,-\, 1)/2,
    &\!\!
    n \in \{1,\ldots,N - 1\},
    \\
    b_N^{(1)} \,=\, 0.
    &
    \ea
    \]

To compute the optimal structured feedback gain for larger values of $\veps$, we use $F(\veps)$ obtained from perturbation analysis to initialize Newton's method, as described in Section~\ref{sec.opt_gain}.

\subsection{Newton's method for larger values of $\veps$}
    \label{sec.opt_gain}

In this section, we employ Newton's method developed in~\cite{farlinjovCDC09} to solve the necessary conditions for optimality~\ref{NC1}-\ref{NC3} as $\veps$ is gradually increased to $1$. Newton's method is an iterative descent algorithm for finding local minima in optimization problems~\cite{boyvan04}. Specifically, given an initial stabilizing structured gain $F^0$, a decreasing sequence of the objective function $\{J(F^i)\}$ is generated by updating $F$ according to
    $
        F^{i+1}
         =
        F^i
         +
        s^i \, \tilde{F}^i.
    $
Here, $\tilde{F}^i$ is the Newton direction that satisfies the structural constraint and $s^i$ is the step-size. The details of computing $\tilde{F}^i$ and choosing the step-size $s^i$ can be found in~\cite{farlinjovCDC09}.

For small $\veps$, we initialize Newton's method using $F(\veps)$ obtained from the perturbation expansion up to the first order in $\veps$,
    $
    F(\veps) = F_0 + \veps F_1.
    $
We then increase $\veps$ slightly and use the optimal structured gain resulting from Newton's method at the previous $\veps$ to initialize the next round of iterations. We continue increasing $\veps$ gradually until desired value $\veps=1$ is reached, that is, until the optimal structured gain $F$ for the desired $Q_d$ is obtained.

Since the homotopy-based Newton's method solves a family of optimization problems parameterized by $\veps$, the optimal feedback gain is a function of $\veps \in [0,1]$. To see the incremental change relative to the spatially uniform gain $F_0$, we consider the difference between the optimal forward gain $f_n(\veps)$ and the uniform gain $f_n(0)=1$,
    \[
        \tilde{f}_n(\veps)
        \,:=\,
        f_n(\veps)
        \,-\,
        f_n(0)
        \,=\,
        f_n(\veps)
        \,-\,
        1.
    \]

\begin{figure}
    \centering
    \subfloat[]
    {
    \includegraphics[width=0.24\textwidth]{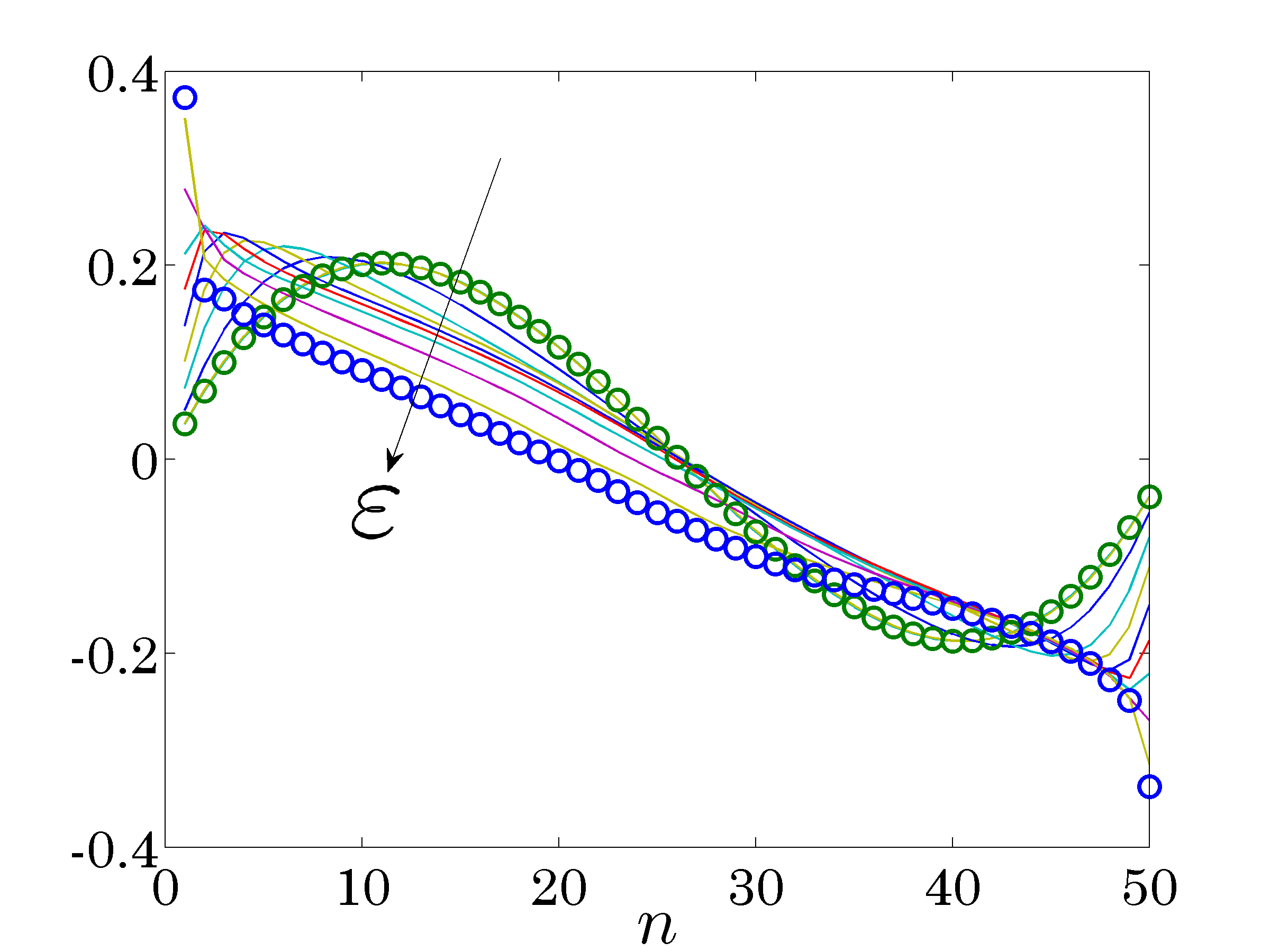}
    \label{fig.forward_gain_Q_I}
    }
    \subfloat[]
    {
    \includegraphics[width=0.24\textwidth]{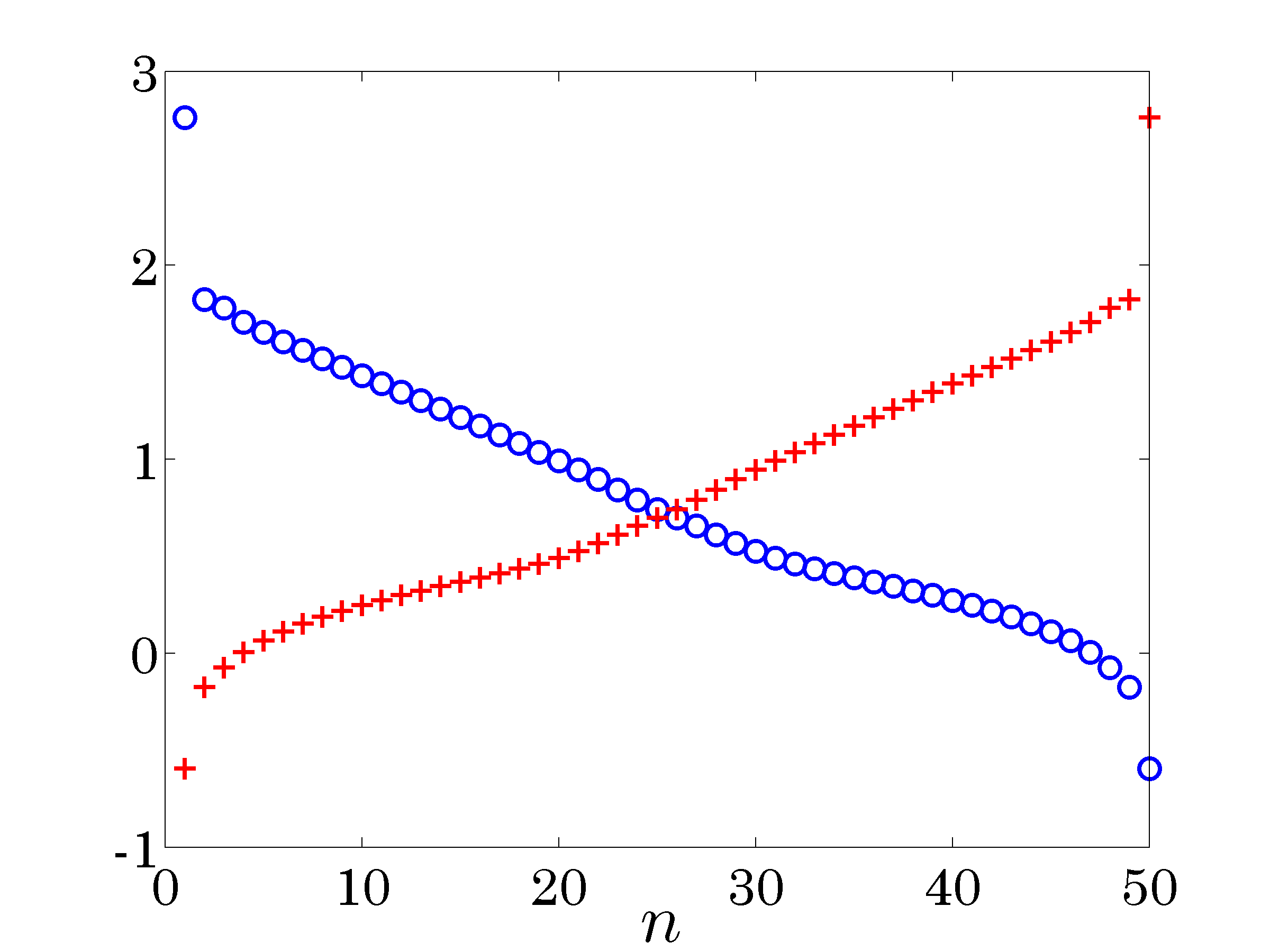}
    \label{fig.fb_gain}
    }
    \caption{Formation with fictitious follower, $N = 50$, $r = 1$, $Q_0 = T^2$, and $Q_d = I$.
    (a) Normalized optimal forward gain $\tilde{f}(\veps) / \| \tilde{f}(\veps) \|$ changes from an almost sinusoidal shape (cf.\ analytical expression in~(\ref{eq.F1})) at $\veps = 10^{-4}$ to an almost piecewise linear shape at $\veps = 1$.
    (b) Optimal forward $(\circ)$ and backward $(+)$ gains at $\veps = 1$.}
    \end{figure}

Figure~\ref{fig.forward_gain_Q_I} shows the normalized profile $\tilde{f}(\veps) / \| \tilde{f}(\veps) \|$ for a formation with fictitious follower, $N = 50$, $r = 1$, $Q_0 = T^2$, and $Q_d = I$. The values of $\veps$ are determined by $20$ logarithmically spaced points between $10^{-4}$ and $1$. As $\veps$ increases, the normalized profile changes from an almost sinusoidal shape (cf.\ analytical expression in~(\ref{eq.F1})) at $\veps = 10^{-4}$ to an almost piecewise linear shape at $\veps = 1$. Note that the homotopy-based Newton's method converges to the same feedback gains at $\veps = 1$ when it is initialized by the optimal symmetric controller obtained in Section~\ref{sec.sym}.

Since the underlying path-graph exhibits symmetry between the edge pairs associated with $f_n$ and $b_{N+1-n}$, the optimal forward and backward gains satisfy a {\em central symmetry\/} property,
    \beq
    \nn
        f_n
        \, = \,
        b_{N+1-n},
        \;\;
        n
        \, \in \,
        \{1,\ldots,N\},
    \eeq
for all $\veps \in [0,1]$; see Fig.~\ref{fig.fb_gain} for $\veps = 1$. We note that the first vehicle has a larger forward gain than other vehicles; this is because it neighbors the fictitious leader. The forward gains decrease as one moves away from the fictitious leader; this is because information about the absolute desired trajectory of the fictitious leader becomes less accurate as it propagates down the formation. Similar interpretation can be given to the optimal backward gains, which monotonically increase as one moves towards the fictitious follower.

\begin{figure}
    \centering
    \subfloat[]
    {
    \includegraphics[width=0.24\textwidth]{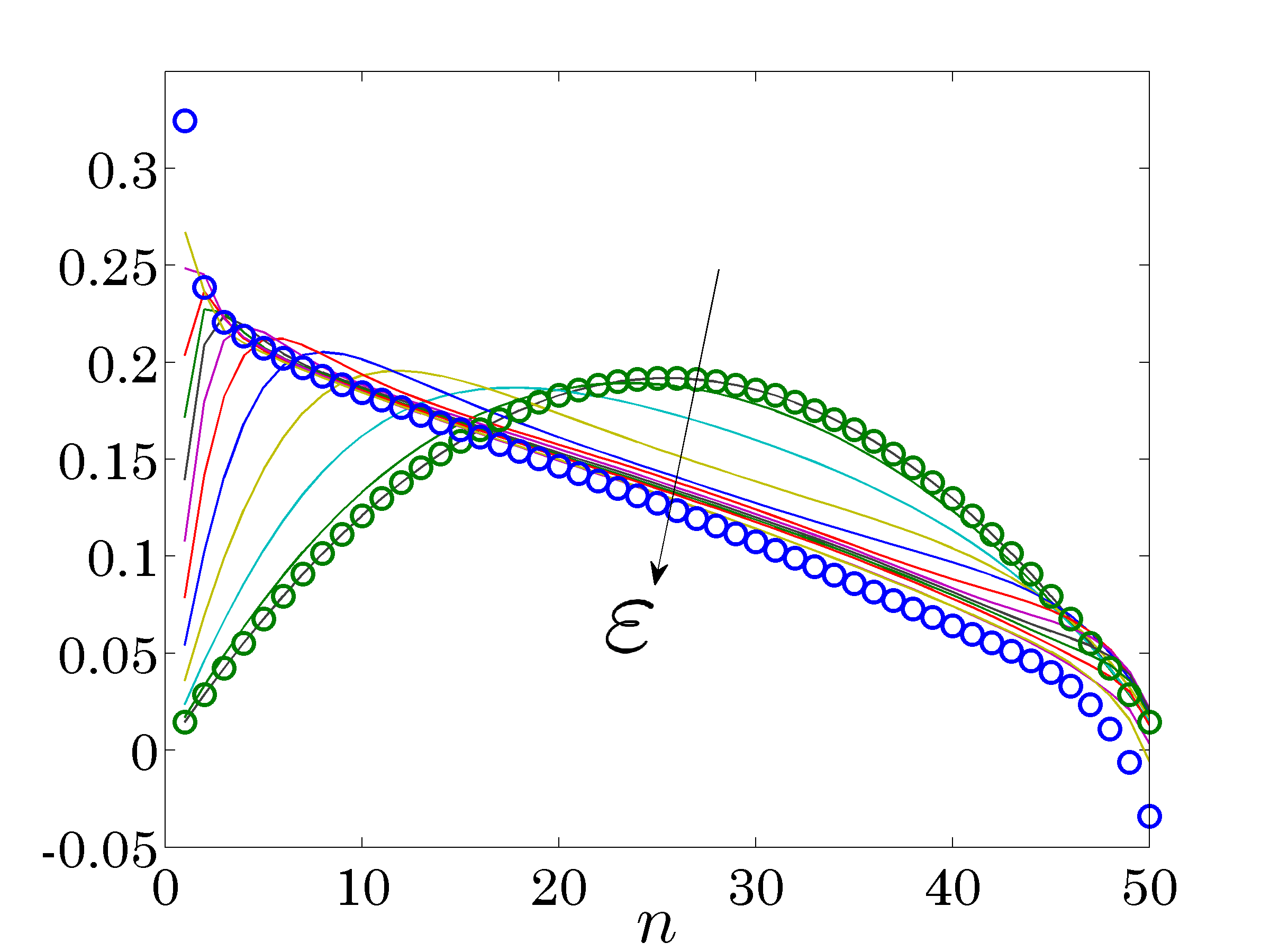}
    \label{fig.forward_gain_Q_I_nofollower}
    }
    \subfloat[]
    {
    \includegraphics[width=0.24\textwidth]{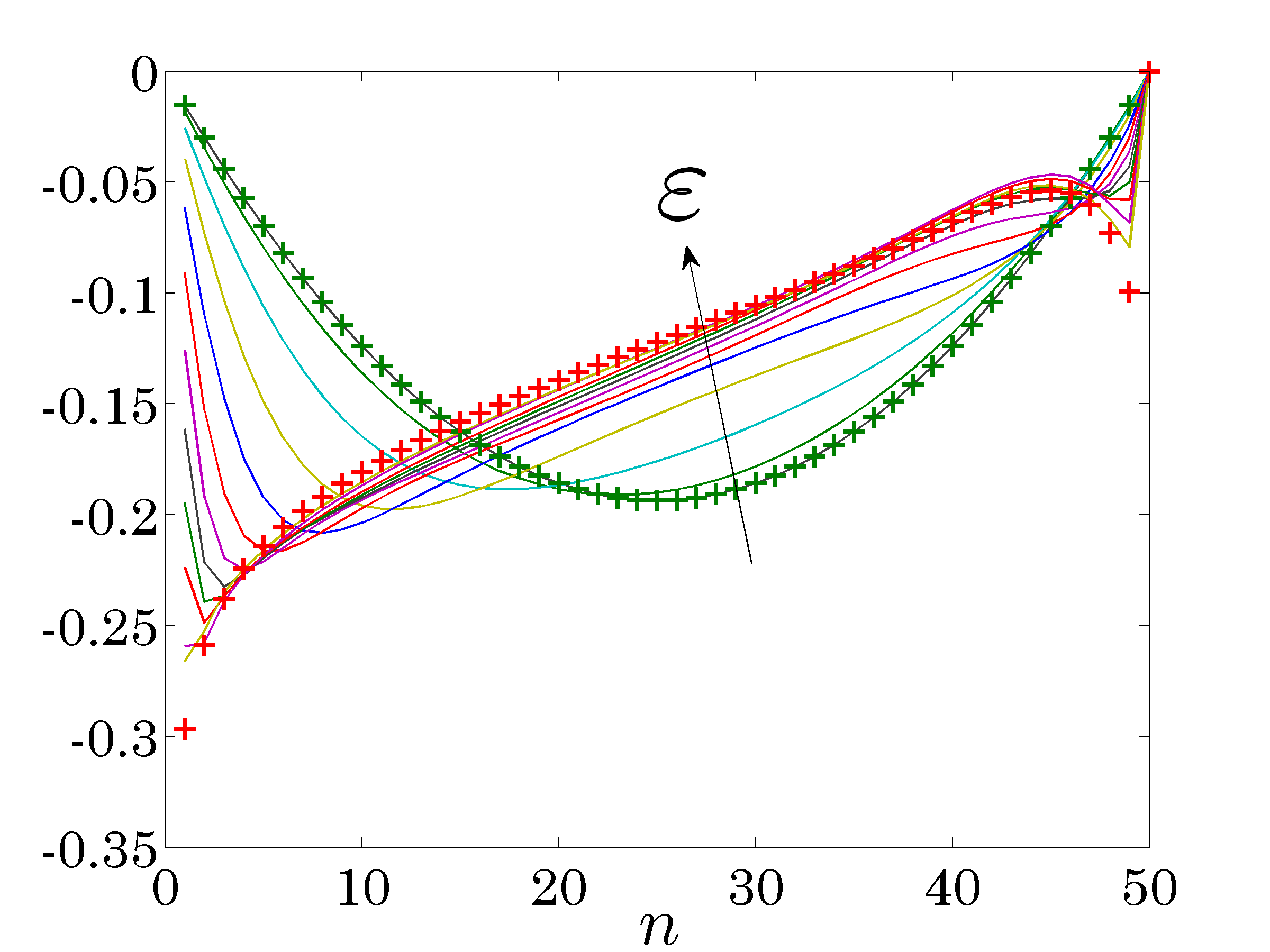}
    \label{fig.backward_gain_Q_I_nofollower}
    }
    \\
    \subfloat[]
    {
    \includegraphics[width=0.24\textwidth]{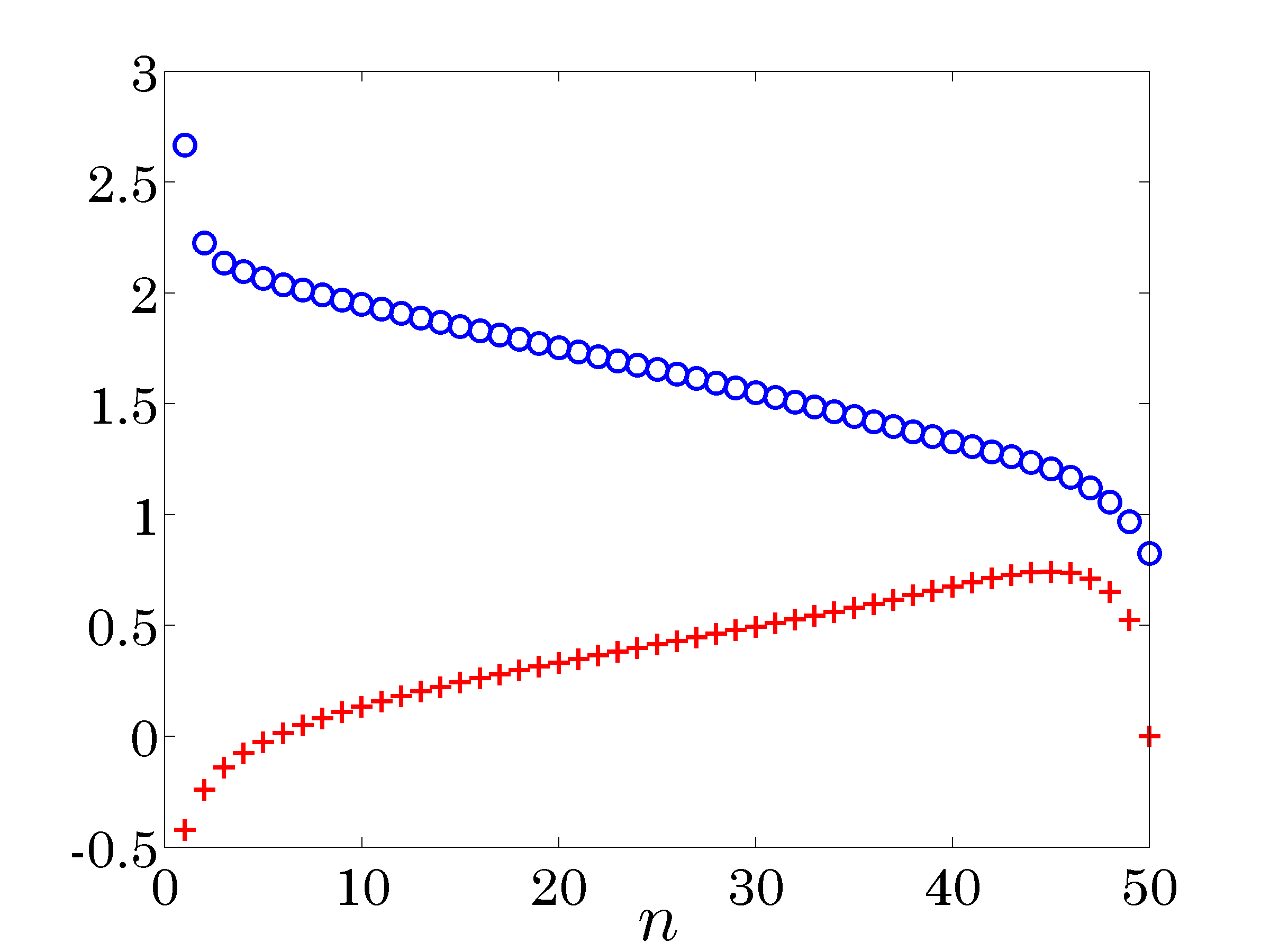}
    \label{fig.feedback_gain_Q_I_nofollower}
    }
    \caption{Formation without fictitious follower, $N = 50$, $r = 1$, $Q_0 = T^2$, and $Q_d = I$.
    Normalized optimal (a) forward and (b) backward gains.
    (c) Optimal forward $(\circ)$ and backward $(+)$ gains at $\veps = 1$.}
    \end{figure}

Since the $1$st vehicle has a {\em negative\/} backward gain (see Fig.~\ref{fig.fb_gain}), if the distance between the $1$st and the $2$nd vehicles is greater than the desired value $\delta$, then the $1$st vehicle distances itself even further from the $2$nd vehicle. On the other hand, if the distance is less than $\delta$, then the $1$st vehicle pulls itself even closer to the $2$nd vehicle. This negative backward gain of the $1$st vehicle can be interpreted as follows: Since the $1$st vehicle has access to its global position, it aims to correct the absolute positions of other vehicles in order to enhance formation coherence. If the $2$nd vehicle is too close to the $1$st vehicle, then the $1$st vehicle moves towards the $2$nd vehicle to push it back; this in turn pushes other vehicles back. If the $2$nd vehicles is too far  from the $1$st vehicle, then the $1$st vehicle moves away from the $2$nd vehicle to pull it forward; this in turn pulls other vehicles forward. Similar interpretation can be given to the negative forward gain of the $N$th vehicle that neighbors the fictitious follower. Also note that the forward gain of the $N$th vehicle becomes {\em positive\/} when the fictitious follower is {\em removed\/} from the formation; see Fig.~\ref{fig.feedback_gain_Q_I_nofollower}. This perhaps suggests that negative feedback gains of the $1$st and the $N$th vehicles are a consequence of the fact that both of them have access to their own global positions.


As shown in Figs.~\ref{fig.forward_gain_Q_I_nofollower} and~\ref{fig.backward_gain_Q_I_nofollower}, the normalized optimal gains for the formation without the fictitious follower also change continuously as $\veps$ increases to $1$. In this case, however, the optimal forward and backward gains do not satisfy the central symmetry; see Fig.~\ref{fig.feedback_gain_Q_I_nofollower}. Since the optimal controller puts more emphasis on the vehicles ahead when the fictitious follower is {\em not\/} present, the forward gains have larger magnitudes than the backward gains. As in the formations with the fictitious follower, the optimal forward gains decrease monotonically as one moves away from the fictitious leader. On the other hand, the optimal backward gains at first increase as one moves away from the $1$st vehicle and then decrease as one approaches the $N$th vehicle in order to satisfy the constraint $b_N = 0$.

\section{Performance vs.\ size for the single-integrator model}
    \label{sec.scaling}

In this section, we study the performance of the optimal symmetric and non-symmetric gains obtained in Sections~\ref{sec.sym} and~\ref{sec.opt_gain}.
This is done by examining the dependence on the formation size of performance measures $\Pi_g$, $\Pi_l$, and $\Pctr$ introduced in Section~\ref{sec.Pi}. Our results highlight the role of non-symmetry and spatial variations on the scaling trends in large-scale formations. They also illustrate performance improvement achieved by the optimal controllers relative to spatially uniform symmetric and non-symmetric feedback gains.

For the spatially uniform {\em symmetric\/} gain with $f_n = b_n = \alpha > 0$, we show {\em analytically\/} that $\Pi_g$ is an {\em affine\/} function of $N$. This implies that the formation coherence scales linearly with $N$ irrespective of the value of $\alpha$. We also {\em analytically\/} establish that the spatially uniform {\em non-symmetric\/} gain with $\{ f_n = \alpha > 0$, $b_n = 0 \}$ (look-ahead strategy) provides a {\em square-root\/} asymptotic dependence of $\Pi_g$ on $N$. Thus, symmetry breaking between the forward and backward gains may improve coherence of large-scale formations. Note that the forward-backward asymmetry also provides more favorable scaling trends of the least damped mode of the closed-loop system~\cite{barmehhes09}. We then investigate how {\em spatially varying\/} optimal feedback gains, introduced in Sections~\ref{sec.sym} and~\ref{sec.opt_gain}, influence coherence of the formation. We show that the optimal symmetric gain provides a {\em square-root\/} dependence of $\Pi_g$ on $N$ and that the optimal non-symmetric gain provides a {\em fourth-root\/} dependence of $\Pi_g$ on $N$.

Even though we are primarily interested in asymptotic scaling of the global performance measure $\Pi_g$, we also examine the local performance measure $\Pi_l$ and the control energy $\Pctr$. From Section~\ref{sec.Pi} we recall that the global and local performance measures quantify the formation-size-normalized ${\cal H}_2$ norm of the transfer function from $d$ to $\zeta_1$ of the closed-loop system,
    \beq
    \ba{rcl}
    \dot{x}
    & \!\! = \!\! &
    - \, F C \, x
    \,+\,
    d
    \\[0.cm]
    \zeta
    & \!\! = \!\! &
    \tbo{\zeta_1}{\zeta_2}
    \,=\,
    \tbo{Q_s^{1/2}}{- \, F C}
    x,
    ~~
    Q_s
    \, = \,
    \left\{
    \ba{rl}
    I, & s \, = \, g,
    \\
    T, & s \, = \, l,
    \ea
    \right.
    \ea
    \non
    \eeq
and that $\Pctr$ is the formation-size-normalized ${\cal H}_2$ norm of the transfer function from $d$ to $\zeta_2$. These can be determined from
    \beq
    \label{eq.Pi_s}
    \ba{rcl}
    \Pi_s
    \!\!& = &\!\!
    (1/N) \, \trace \left( L \, Q_s \right),
    \\[0.1cm]
    \Pctr
    \!\!& = &\!\!
    (1/N) \, \trace \left( L \, C^T F^T F C \right),
    \ea
    \eeq
where $L$ denotes the closed-loop controllability Gramian,
    \beq
    \left( - \, F C \right) L
    \,+\,
    L \left( - \, F C \right)^T
    \,=\,
    - \, I.
    \label{eq.Lctrb}
    \eeq

The asymptotic scaling properties of $\Pi_g$, $\Pi_l$, and $\Pctr$, for the above mentioned spatially uniform controllers and the spatially varying optimal controllers, obtained by solving~\ref{SH2} with $Q = I$ and $r = 1$,
are summarized in Table~\ref{tab.scaling}. For both spatially uniform symmetric and look-ahead strategies, we {\em analytically\/} determine the dependence of these performance measures on the formation size in Sections~\ref{sec.uni_sym_ctr} and~\ref{sec.uni_nonsym_ctr}. Furthermore, for the formation without the fictitious follower subject to the optimal symmetric gains, we provide {\em analytical\/} results in Section~\ref{sec.opt_ctr}. For the optimal symmetric and non-symmetric gains in the presence of fictitious followers, the scaling trends are obtained with the aid of numerical computations in Section~\ref{sec.opt_ctr}.

Several comments about the results in Table~\ref{tab.scaling} are given next. First, in contrast to the spatially uniform controllers, the optimal symmetric and non-symmetric gains, resulting from an $N$-independent control penalty $r$ in~\ref{SH2}, do not provide uniform bounds on the control energy per vehicle, $\Pctr$. This implies the trade-off between the formation coherence $\Pi_g$ and control energy $\Pctr$ in the design of the optimal controllers. It is thus of interest to examine formation coherence for optimal controllers with bounded control energy per vehicle (see Remark~\ref{rem.bdd_ctr}). Second, the controller structure (e.g., symmetric or non-symmetric gains) plays an important role in the formation coherence. In particular, departure from symmetry in localized feedback gains can significantly improve coherence of large-scale formations (see Remark~\ref{rem.ctr1}).

    \begin{table*}
\caption{Asymptotic dependence of $\Pi_g$, $\Pi_l$, and $\Pctr$ on the formation size $N$ for uniform symmetric, uniform non-symmetric (look-ahead strategy), and optimal symmetric and non-symmetric gains of Sections~\ref{sec.sym} and~\ref{sec.opt_gain} with $Q = I$ and $r = 1$. The scalings displayed in red are determined {\em analytically\/}; other scalings are estimated based on numerical computations.}
  \label{tab.scaling}
  \centering
    \begin{tabular}{lccc}
         {\bf Controller}                              & $\Pi_g$                       & $\Pi_l$                        &  $\Pctr$      \\[0.1cm]
         uniform symmetric with/without follower       & \tc{red}{$O(N)$}             & \tc{red}{$O(1)$}              &  \tc{red}{$O(1)$}               \\[0.1cm]
         uniform non-symmetric                         & \tc{red}{$O(\sqrt{N})$}      & \tc{red}{$O(1)$}              &  \tc{red}{$O(1)$}               \\[0.1cm]
         optimal symmetric without follower            & \tc{red}{$O(\sqrt{N})$}      & \tc{red}{$O(1/\sqrt{N})$}     &  \tc{red}{$O(\sqrt{N})$}        \\[0.1cm]
         optimal symmetric with follower               & $O(\sqrt{N})$      & $O(1/\sqrt{N})$     &  $O(\sqrt{N})$        \\[0.1cm]
         optimal non-symmetric with/without follower   & $O(\sqrt[4]{N})$   & $O(1/\sqrt[4]{N})$  &  $O(\sqrt[4]{N})$
    \end{tabular}
    \end{table*}

\subsection{Spatially uniform symmetric gain}
    \label{sec.uni_sym_ctr}

For the spatially uniform symmetric controller with $f_n = b_n = \alpha > 0$, we next show that $\Pi_g$ is an affine function of $N$ and that, in the limit of an infinite number of vehicles, both $\Pi_l$ and $\Pctr$ become formation-size-independent. These results hold irrespective of the presence of the fictitious follower.

For the single-integrator model with the fictitious follower we have $K = F C = \alpha T$ (see (\ref{eq.T}) for the definition of $T$), and
    $
    L
    =
    T^{-1}/(2 \alpha )
    $
solves the Lyapunov equation~(\ref{eq.Lctrb})~\cite[Lemma~1]{bamdah03}. Since the $n$th diagonal entry of $T^{-1}$ is determined by~(cf.~(\ref{eq.invK}))
    \[
    (T^{-1})_{nn}
     \,=\,
    n \left( N + 1 - n \right)/ \left( N+1 \right),
    \]
from~(\ref{eq.Pi_s}) we conclude that the global performance measure $\Pi_g$ is an affine function of $N$, and that both $\Pi_l$ and $\Pctr$ are formation-size-independent,
    \beq
    \ba{rcl}
    \Pi_g
    \!\! & = & \!\!
    \trace
    \left(
    T^{-1}
    \right)
    /(2 \alpha N)
    \\
    \!\! & = & \!\!
    \dfrac{1}{2 \alpha N} \ds{\sum_{n \, = \, 1}^{N}} \, n
    \, - \,
    \dfrac{1}{2 \alpha N (N+1)} \ds{\sum_{n \, = \, 1}^{N}} \, n^2
    \, = \,
    \dfrac{N + 2}{12 \alpha },
    \\[0.45cm]
    \Pi_l
    \!\! & = & \!\!
    \trace
    \left(
    T \, T^{-1}
    \right)/(2 \alpha N)
    \, = \,
    1/ (2 \alpha),
    \\[0.1cm]
    \Pctr
    \!\! & = & \!\!
    \trace
    \left(
    \alpha^2 \, T \, T \, T^{-1}
    \right)
    /
    (2 \alpha N)
    \,=\,
    \alpha.
    \ea
    \non
    \eeq
For the formation without the fictitious follower, the following expressions
    \[
    \Pi_g
    \,=\,
    (N+1)/(4\alpha),
    ~~
    \Pi_l
    \,=\,
    1/\alpha,
    ~~
    \Pctr
    \,=\,
    \alpha(3N+1)/(2N),
    \]
imply that, for the spatially uniform symmetric controller, the asymptotic scaling trends do not depend on the presence of the fictitious follower (derivations omitted for brevity).

\subsection{Spatially uniform non-symmetric gain~(look-ahead strategy)}
    \label{sec.uni_nonsym_ctr}

We next examine the asymptotic scaling of the performance measures for the spatially uniform non-symmetric gain with $\{ f_n = \alpha > 0$, $b_n = 0 \}$. We establish the {\em square-root\/} scaling of $\Pi_g$ with $N$ and the formation-size-independent scaling of $\Pi_l$. Furthermore, in the limit of an infinite number of vehicles, we show that $\Pctr$ becomes $N$-independent.

For the single-integrator model with $K = F C = \alpha C_f$ (see (\ref{eq.Cf}) for the definition of $C_f$), the solution of the Lyapunov equation~(\ref{eq.Lctrb}) is given by
    \beq
    L
    \,=\,
    \int_0^\infty
    \mre ^{ - \,  \alpha \, C_f \, t }
    \,
    \mre ^{ - \,  \alpha \, C_f^T \, t }
    \, \mrd t.
    \label{eq.Lint}
    \eeq
As shown in Appendix~\ref{app.nonsymmetric}, the inverse Laplace transform of \mbox{$(sI + \alpha C_f)^{-1}$} can be used to determine the analytical expression for $\mre^{ -\,  \alpha \, C_f \, t }$, yielding the following formulae,
    \beq
    \non
    \ba{rcl}
    \Pi_g (N)
    \!\! & = & \!\!
    \dfrac{1}{N}
    \ds{\sum_{n\,=\,1}^N}
    L_{nn}
    \, = \,
    \dfrac{1}{N}
    \ds{\sum_{n\,=\,1}^N}
    \dfrac{\alpha \, \Gamma(n + 1/2)}{\sqrt{\pi} \, \Gamma(n)}
    \\[0.5cm]
    \!\! & = & \!\!
    \dfrac{2 \, \alpha \, \Gamma(N + 3/2)}{3 \, \sqrt{\pi} \, \Gamma(N + 1)},
    \\
    \Pi_l
    \!\! & = & \!\!
    \alpha,
    \\
    \Pctr
    \!\! & = & \!\!
    \alpha
    \, - \,
    (1/N) L_{NN},
    \ea
    \eeq
with $\Gamma(\cdot)$ denoting the Gamma function. These are used in Appendix~\ref{app.nonsymmetric} to show that, in the limit of an infinite number of vehicles, a look-ahead strategy for the single-integrator model provides the square-root dependence of $\Pi_g$ on $N$ and the formation-size-independent $\Pi_l$ and $\Pctr$.

\subsection{Optimal symmetric and non-symmetric controllers}
    \label{sec.opt_ctr}

We next examine the asymptotic scaling of the performance measures for the optimal symmetric and non-symmetric gains of Sections~\ref{sec.sym} and~\ref{sec.opt_gain}. For the formation without the fictitious follower, we {\em analytically\/} establish that the optimal symmetric gains asymptotically provide $O(\sqrt{N})$, $O(1/\sqrt{N})$, and $O(\sqrt{N})$ scalings of $\Pi_g$, $\Pi_l$, and $\Pctr$, respectively. We then use numerical computations to (i) confirm these scaling trends for the optimal symmetric gains in the presence of the fictitious follower; and to (ii) show a fourth-root dependence of $\Pi_g$ and $\Pctr$ on $N$ and an $O (1/\sqrt[4]{N})$ dependence of $\Pi_l$ for the optimal non-symmetric gains. All these scalings are obtained by solving~\ref{SH2} with the formation-size-independent control penalty $r$ and $Q = I$. We also demonstrate that uniform control variance (per vehicle) can be obtained by judicious selection of an $N$-dependent $r$. For the optimal symmetric and non-symmetric gains, this constraint on control energy (variance) increases the asymptotic dependence of $\Pi_g$ on $N$ to linear and square-root, respectively.

For the formation without the fictitious follower, the optimal symmetric gains are given by~(\ref{eq.symK_QI}). As shown in~(\ref{eq.KequKinv}), $\trace \, (K^{-1}) = \trace \, ( r K)$, thereby yielding
    \beq
    \ba{rcl}
        \Pi_g
        \!\! & = & \!\!
        r
        \,
        \Pctr
        \,=\, \dfrac{1}{2N}
        \,
        \trace \left( K^{-1} \right)
        \\
        \!\! & = & \!\!
        \dfrac{\sqrt{r}}{2N}
        \left( \sqrt{N}
        \, + \,
        \ds{\sum_{n\,=\,1}^{N-1}} \sqrt{2n} \right).
        \ea
    \label{eq.Pg=Pctr}
    \eeq
In the limit of an infinite number of vehicles,
    \[
    \ba{rcl}
    \ds{\lim_{N \, \to \, \infty}}
        \dfrac{\Pi_g (N)}{\sqrt{N}}
    \!\! & = & \!\!
    \ds{\lim_{N \, \to \, \infty}}
        \ds{\sum_{n\,=\,1}^{N-1}} \sqrt{\dfrac{r n}{2N}}
        \,
        \dfrac{1}{N}
        \\[0.45cm]
    \!\! & = & \!\!
        \ds{\int_{0}^{1}}  \sqrt{\dfrac{r x}{2}}  \, \mrd x
        \,=\,
        \sqrt{\dfrac{2 r}{9}},
    \ea
    \]
which, for an $N$-independent $r$, leads to an asymptotic square-root dependence of $\Pi_g$ and $\Pctr$ on $N$,
    \beq
    \label{eq.pigpctr}
    \ba{rcl}
    \Pi_g (N)
    \!\! & = & \!\!
    \sqrt{\dfrac{2 \, r N}{9}}
    \, + \,
    \sqrt{\dfrac{r}{4 N}},
    ~~
    N
    \, \gg \,
    1
    \\
    \Pctr (N)
    \!\! & = & \!\!
    \sqrt{\dfrac{2N}{9r}}
    \, + \,
    \dfrac{1}{\sqrt{4 \, r N}},
    ~~
    N
    \, \gg \,
    1.
    \ea
    \eeq
Similar calculation can be used to obtain $O(1/\sqrt{N})$ asymptotic scaling of $\Pi_l$.

\begin{figure}
      \centering
        \subfloat[]
        {\label{fig.rootfit_Pi_g}
        \includegraphics[width=0.24\textwidth]{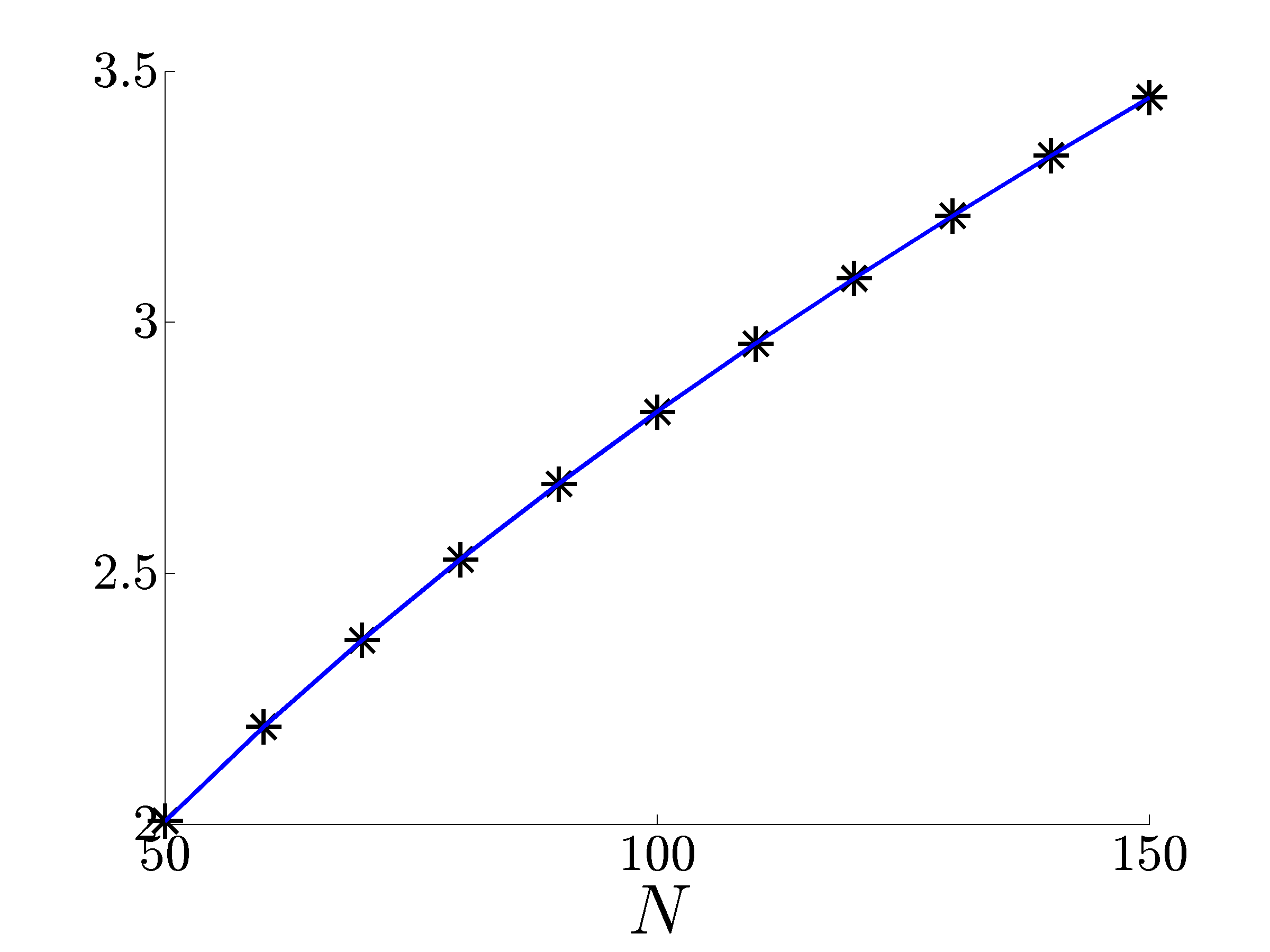}}
        \subfloat[]
        {\label{fig.root4fit_Pi_g}
        \includegraphics[width=0.24\textwidth]{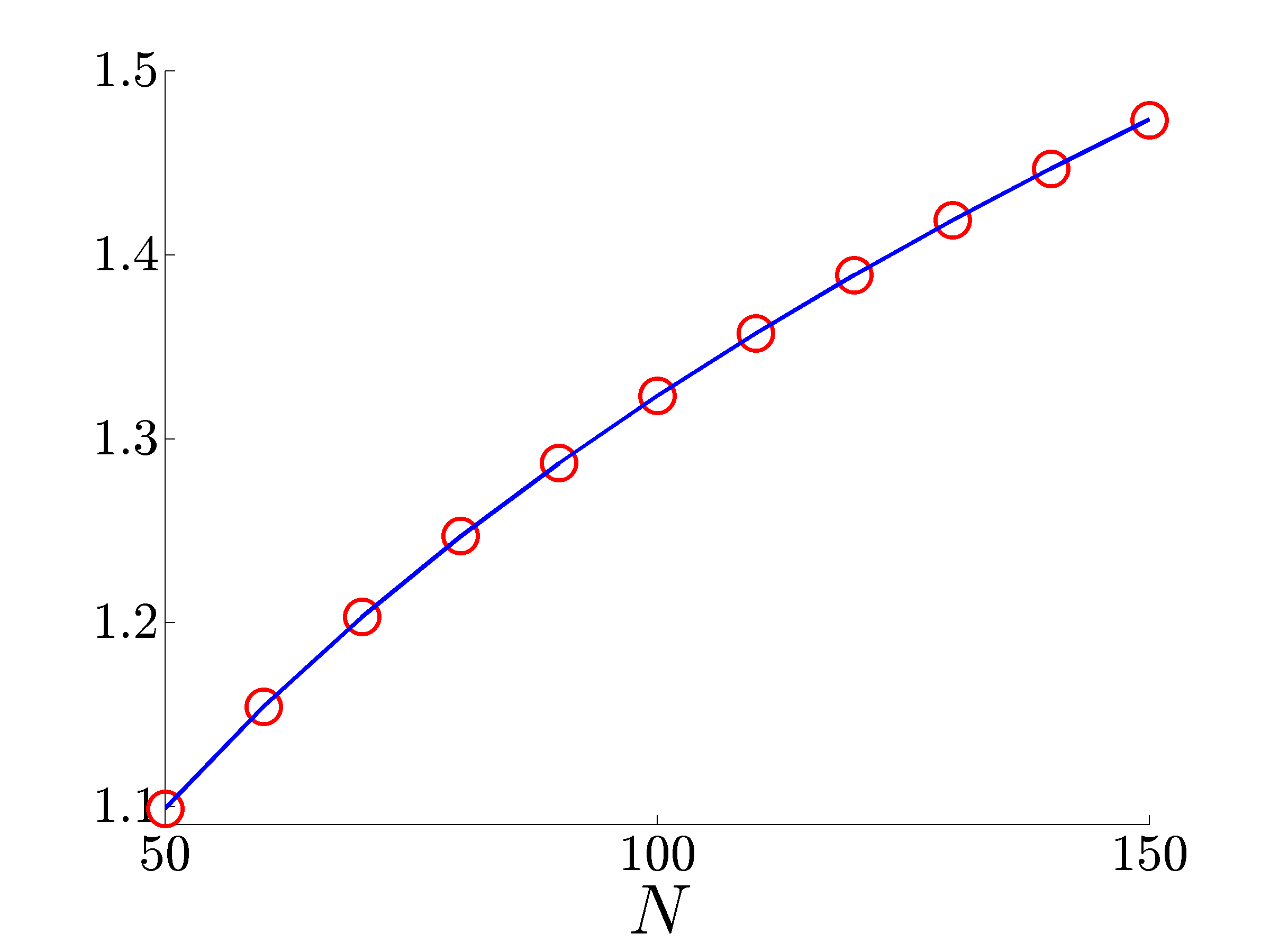}}
      \caption{(a) Square-root scaling of $\Pi_g \, (*)$ using optimal symmetric gain of Section~\ref{sec.sym}, $0.2784 \sqrt{N}+0.0375$ (curve); and
      (b) Fourth-root scaling of $\Pi_g \, (\circ)$ using optimal non-symmetric gain of Section~\ref{sec.opt_gain}, $0.4459\sqrt[4]{N}-0.0866$ (curve).
      The optimal controllers are obtained by solving~\ref{SH2} with $Q = I$ and $r = 1$ for the formation with the fictitious follower.}
      \label{fig.sqrtfit_Pi_g}
    \end{figure}

\begin{figure}
      \centering
        \subfloat[]
        {\label{fig.rootfit_Pi_l}
        \includegraphics[width=0.24\textwidth]{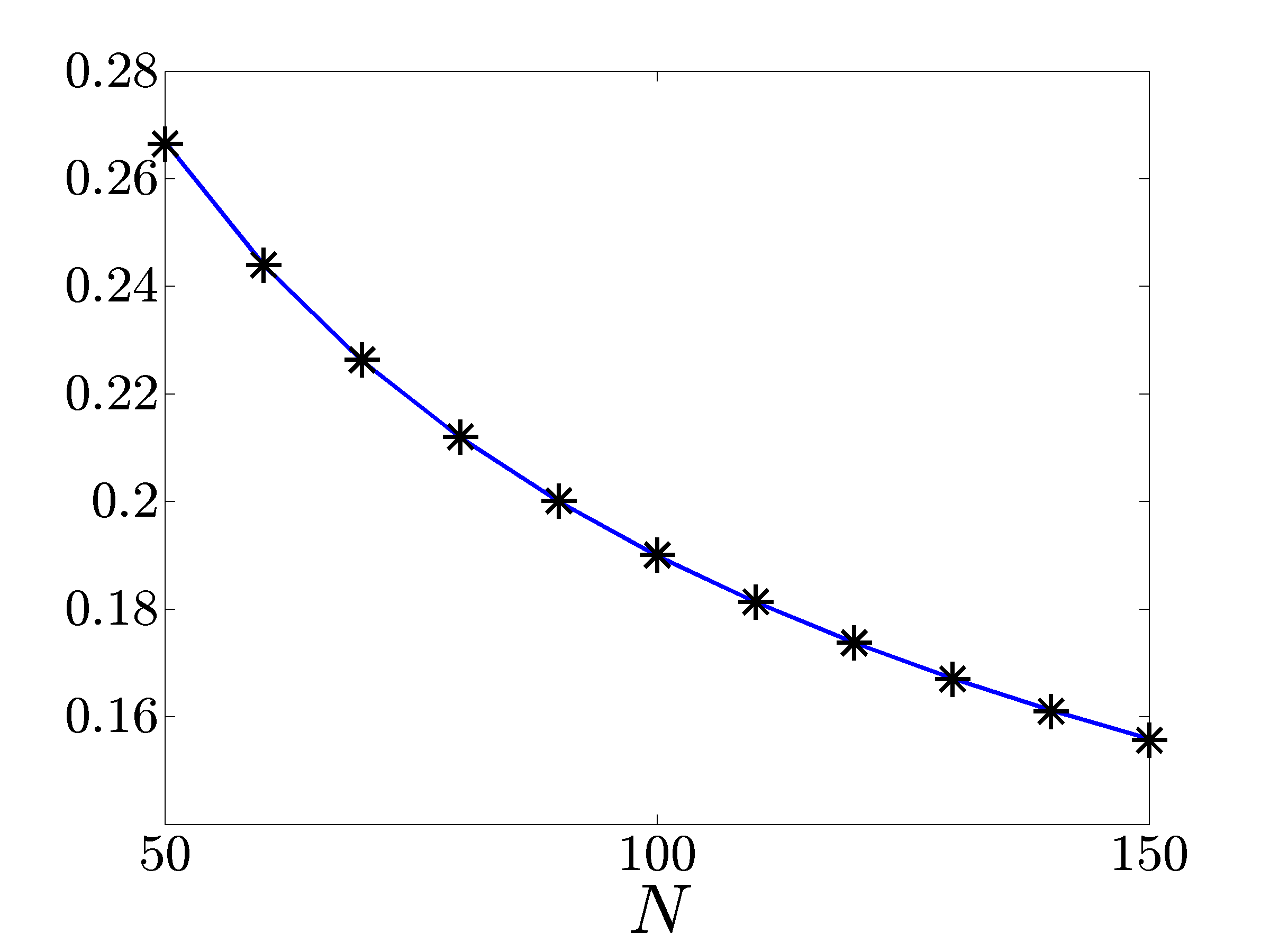}}
        \subfloat[]
        {\label{fig.root4fit_Pi_l}
        \includegraphics[width=0.24\textwidth]{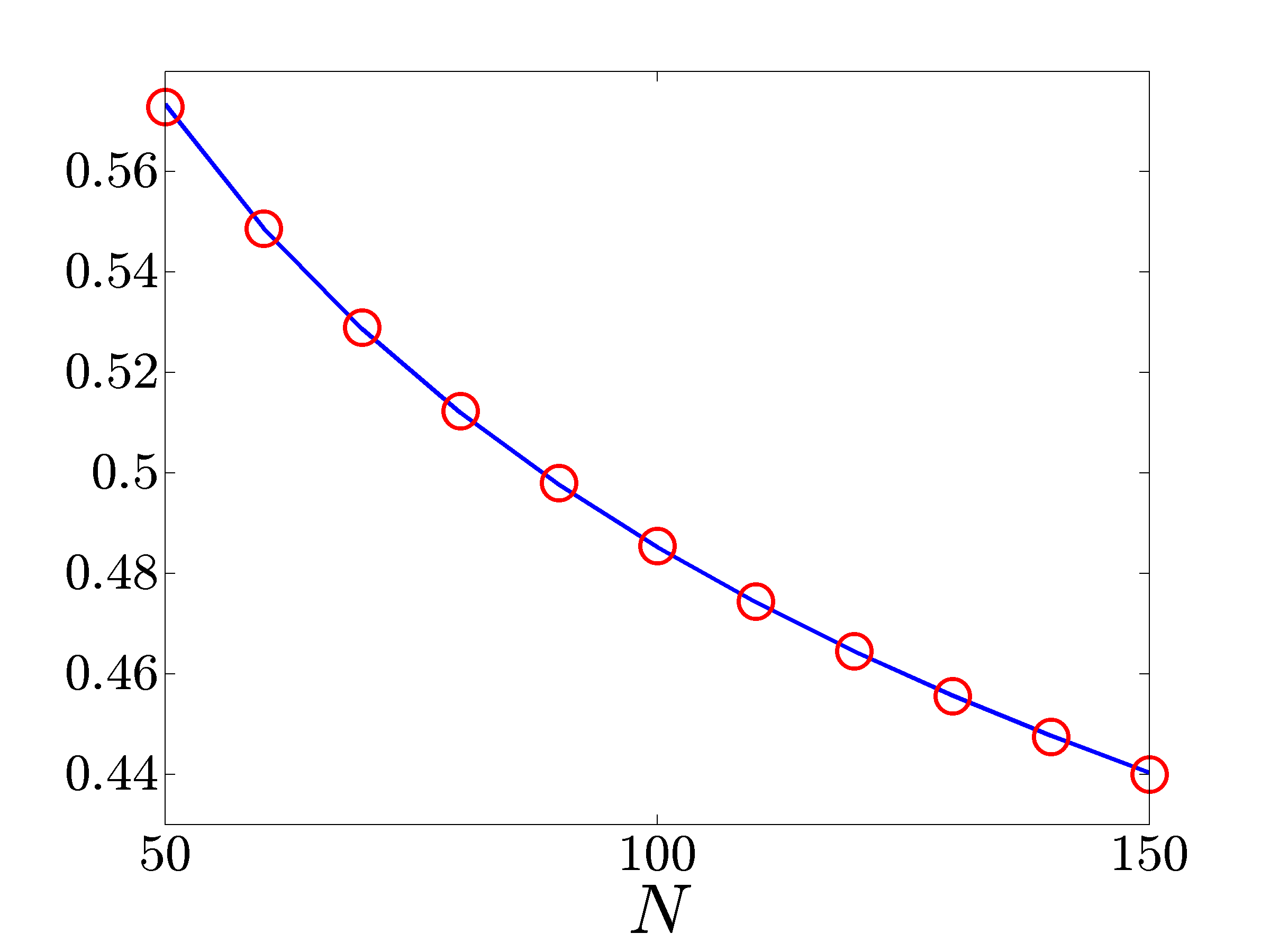}}
      \caption{(a) $\Pi_l \, (*)$ using the optimal symmetric gain of Section~\ref{sec.sym}, $1.8570/\sqrt{N} + 0.0042$ (curve);
      and
               (b) $\Pi_l \, (\circ)$ using the optimal non-symmetric gain of Section~\ref{sec.opt_gain}, $1.4738/\sqrt[4]{N} + 0.0191$ (curve). The optimal controllers are obtained by solving~\ref{SH2} with $Q = I$ and $r = 1$ for the formation with the fictitious follower.}
      \label{fig.sqrtfit_Pi_l}
    \end{figure}

We next use numerical computations to study the scaling trends for the optimal symmetric and non-symmetric gains in the presence of fictitious followers. The optimal symmetric gain (cf.\ ($\circ$) in Fig.~\ref{fig.sym_gain}) provides a {\em square-root\/} scaling of $\Pi_g$ with $N$; see Fig.~\ref{fig.rootfit_Pi_g}. On the other hand, the optimal non-symmetric gain (cf.\ Fig.~\ref{fig.fb_gain}) leads to a {\em fourth-root\/} scaling of $\Pi_g$ with $N$; see Fig.~\ref{fig.root4fit_Pi_g}. The local performance measure $\Pi_l$ decreases monotonically with $N$ for both controllers, with $\Pi_l$ scaling as $O(1/\sqrt{N})$ for the optimal symmetric gain and as $O(1/\sqrt[4]{N})$ for the optimal non-symmetric gain; see Fig.~\ref{fig.sqrtfit_Pi_l}. For both the optimal symmetric and non-symmetric controllers, our computations indicate equivalence between the control energy and the global performance measure when $r = 1$.
(For the optimal symmetric gain without the fictitious follower and $r = 1$, we have {\em analytically\/} shown that
    $
    \Pctr
    =
    \Pi_g;
    $
see formula~(\ref{eq.Pg=Pctr}).) Therefore, the asymptotic scaling of the formation-size-normalized control energy is $O(\sqrt{N})$ for the optimal symmetric gain and $O(\sqrt[4]{N})$ for the optimal non-symmetric gain. Finally, for the formations without the fictitious follower, our computations indicate that the optimal non-symmetric gains also asymptotically provide $O(\sqrt[4]{N})$, $O(1/\sqrt[4]{N})$, and $O(\sqrt[4]{N})$ scalings of $\Pi_g$, $\Pi_l$, and $\Pctr$, respectively.

\begin{figure}
    \begin{center}
    \includegraphics[width=0.45\textwidth]{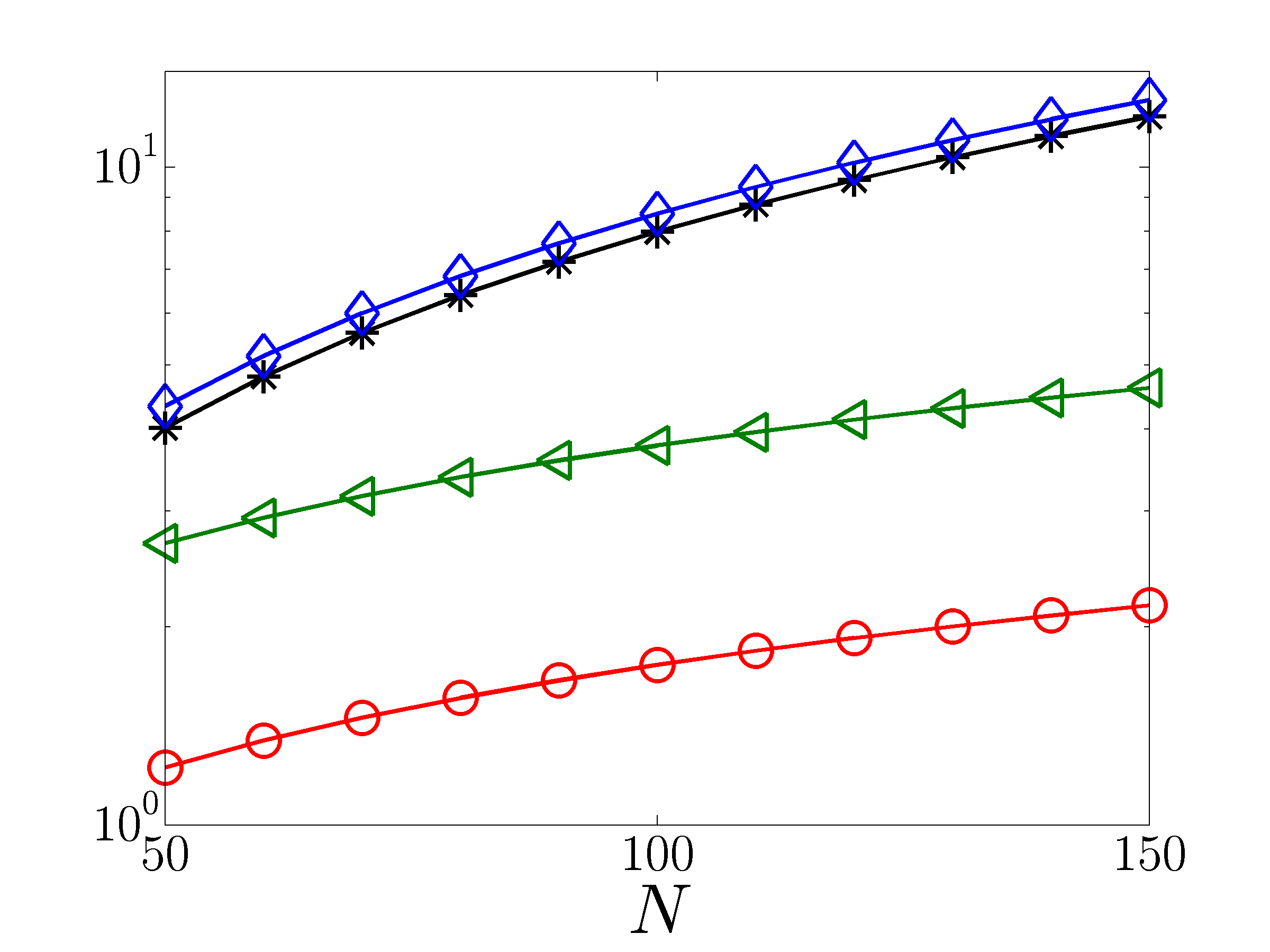}
    \caption{$\Pi_g$ using four structured gains with $\Pctr \approx 1$ for formations with fictitious follower: spatially uniform symmetric ($\diamond$), $N/12 + 1/6$ (blue curve), spatially uniform non-symmetric ($\triangleleft$), $2\sqrt{N} /(3 \sqrt{\pi}) $ (green curve), optimal symmetric ($*$), $0.0793 N + 0.0493$ (black curve), and optimal non-symmetric ($\circ$), $0.1807 \sqrt{N} - 0.0556$ (red curve).}
    \label{fig.Pig_ctr1}
    \end{center}
    \end{figure}

    \remark
    \label{rem.bdd_ctr}
In contrast to the spatially uniform controllers, the optimal structured controllers of Sections~\ref{sec.sym} and~\ref{sec.opt_gain}, resulting from an $N$-independent control penalty $r$ in~\ref{SH2}, do not provide uniform bounds on the formation-size-normalized control energy. These controllers are obtained using ${\cal H}_2$ framework in which control effort represents a `soft constraint'. It is thus of interest to examine formation coherence for optimal controllers with bounded control energy per vehicle. For formations without the fictitious follower, from~(\ref{eq.pigpctr}) we see that the optimal symmetric controller with $r(N) = 2N/9$ asymptotically yields $\Pctr \approx 1$ and $\Pi_g \approx 2N/9 \, + \, 1/(3 \sqrt{2})$. Similarly, for formations with followers, the optimal gains that result in $\Pctr \approx 1$ for large $N$ can be obtained by changing control penalty from $r = 1$ to $r(N) = 0.08 N$ for the optimal symmetric gain and to $r(N) = 0.175 \sqrt{N}$ for the optimal non-symmetric gain\footnote{Both spatially uniform symmetric and look-ahead strategies with $\alpha = 1$ yield $\Pctr = 1$ in the limit of an infinite number of vehicles.}. These $N$-dependent control penalties provide an {\em affine\/} scaling of $\Pi_{g}$ with $N$ for the optimal symmetric gain and a {\em square-root\/} scaling of $\Pi_g$ with $N$ for the optimal non-symmetric gain; see Fig.~\ref{fig.Pig_ctr1}. The asymptotic scalings for formations without followers subject to the optimal symmetric gains are obtained analytically (cf.\ (\ref{eq.pigpctr})); all other scalings are obtained with the aid of computations.

    \remark
    \label{rem.ctr1}
Figure~\ref{fig.Pig_ctr1} illustrates the global performance measure $\Pi_g$ obtained with four aforementioned structured controllers that asymptotically yield $\Pctr \approx 1$ for formations with fictitious follower. Note that the simple look-ahead strategy outperforms the optimal symmetric gain; $O(\sqrt{N})$ vs.\ $O(N)$ scaling. Thus, departure from symmetry in localized feedback gains can significantly improve coherence of large-scale formations. In particular, we have provided an example of a spatially uniform non-symmetric controller that yields better scaling trends than the optimal spatially varying controller obtained by restricting design to symmetric gains. Given the extra degrees of freedom in the optimal symmetric gain this is perhaps a surprising observation, indicating that the network topology may play a more important role than the optimal selection of the feedback gains in performance of large-scale interconnected systems. On the other hand, our results show that the optimal localized controller that achieves the best performance is both {\em non-symmetric and spatially-varying\/}.

\section{Double-integrator model}
    \label{sec.double}

In this section, we solve~\ref{SH2} for the double-integrator model using the homotopy-based Newton's method. We then discuss the influence of the optimal structured gain on the asymptotic scaling of the performance measures introduced in Section~\ref{sec.Pi}. For a formation in which each vehicle -- in addition to relative positions with respect to its immediate neighbors -- has access to {\em its own velocity\/}, our results highlight similarity between optimal forward and backward position gains for the single- and the double-integrator models. We further show that the performance measures exhibit similar scaling properties to those found in single-integrators. We also establish convexity of~\ref{SH2} for the double-integrator model by restricting the controller to symmetric position and uniform diagonal velocity gains.

The perturbation analysis and the homotopy-based Newton's method closely follow the procedure described in Sections~\ref{sec.pert_ana} and~\ref{sec.opt_gain}, respectively. In particular,
    $
    F_0
    =
    [ \, \alpha I \;\; \alpha I \;\; \beta I \, ]
    $
yields
    $
    K_0 = F_0 \, C = [ \, \alpha T \;\; \beta I \, ].
    $
As shown in~\cite{jovIO10}, for positive $\alpha$ and $\beta$ with $\beta^2 > 8 \alpha$, this spatially uniform structured feedback gain is stabilizing and inversely optimal with respect to
    \[
    \ba{rcl}
    Q_0
    \!\! & = & \!\!
    \tbt{Q_p}{O}{O}{Q_v},
    ~
    Q_p
    \, = \,
    r \alpha^2 \, T^2,
    \\[0.35cm]
    Q_v
    \!\! & = & \!\!
    r (\beta^2 \, I - 2 \, \alpha \, T),
    ~
    r
    \, > \,
    0.
    \ea
    \]

\begin{figure}
    \centering
    \subfloat[]
    {\label{fig.fb_gain_2nd}
    \includegraphics[width=0.24\textwidth]{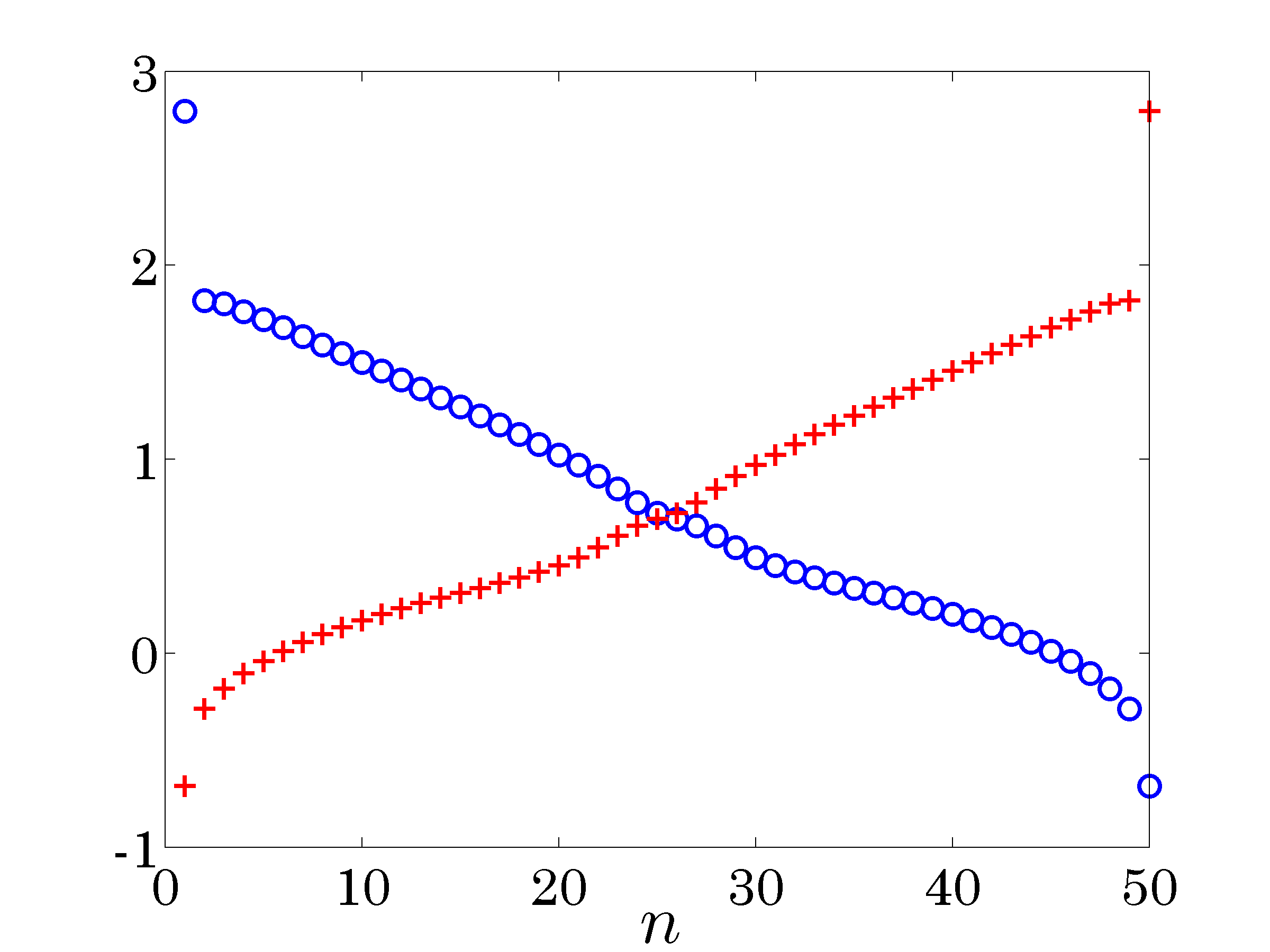}
    }
    \subfloat[]
    {\label{fig.2ndOrder_vel_gain}
    \includegraphics[width=0.24\textwidth]{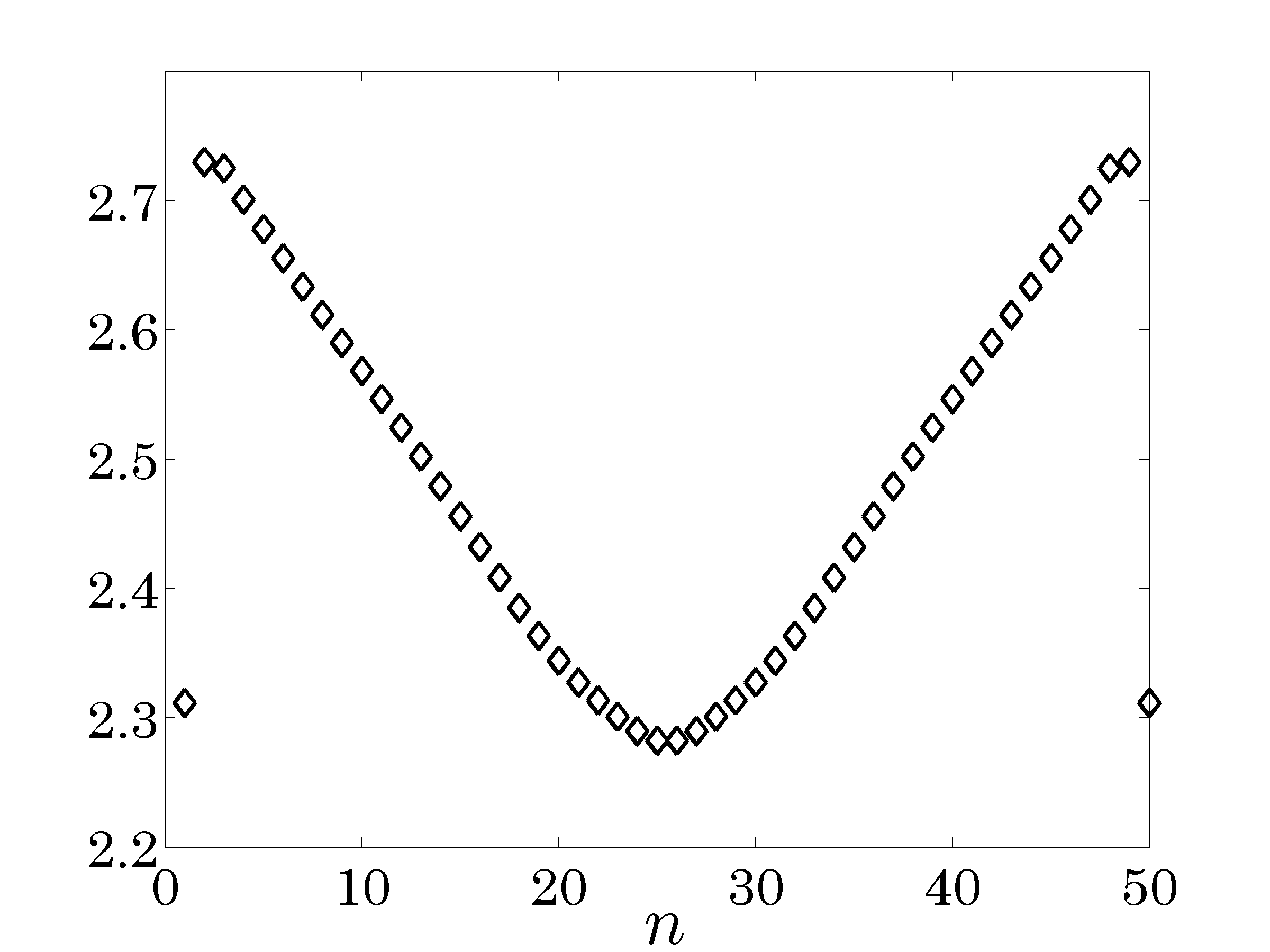}
    }
    \caption{Double-integrator model with fictitious follower, $N = 50$, $Q = I$ and $r = 1$. (a) The optimal forward~($\circ$) and backward gains~($+$); (b) the optimal velocity gains~$(\diamond)$.}
   \end{figure}

In what follows, we choose $\alpha = 1$ and $\beta = 3$ and employ the homotopy-based Newton's method to solve~\ref{SH2} for the double-integrator model.  For a formation with fictitious follower, $N = 50$, $Q = I$, and $r = 1$ the optimal forward and backward position gains are shown in Fig.~\ref{fig.fb_gain_2nd} and the optimal velocity gains are shown in Fig.~\ref{fig.2ndOrder_vel_gain}. We note remarkable similarity between the optimal position gains for the single- and the double-integrator models; cf.\ Fig.~\ref{fig.fb_gain_2nd} and Fig.~\ref{fig.fb_gain}. For a formation without fictitious follower, the close resemblance between the optimal position gains for both models is also observed.

\begin{figure}
    \centering
        \subfloat[$\Pi_g$]
        {\label{fig.Pi_g_Fg_2nd}
        \includegraphics[width=0.24\textwidth]{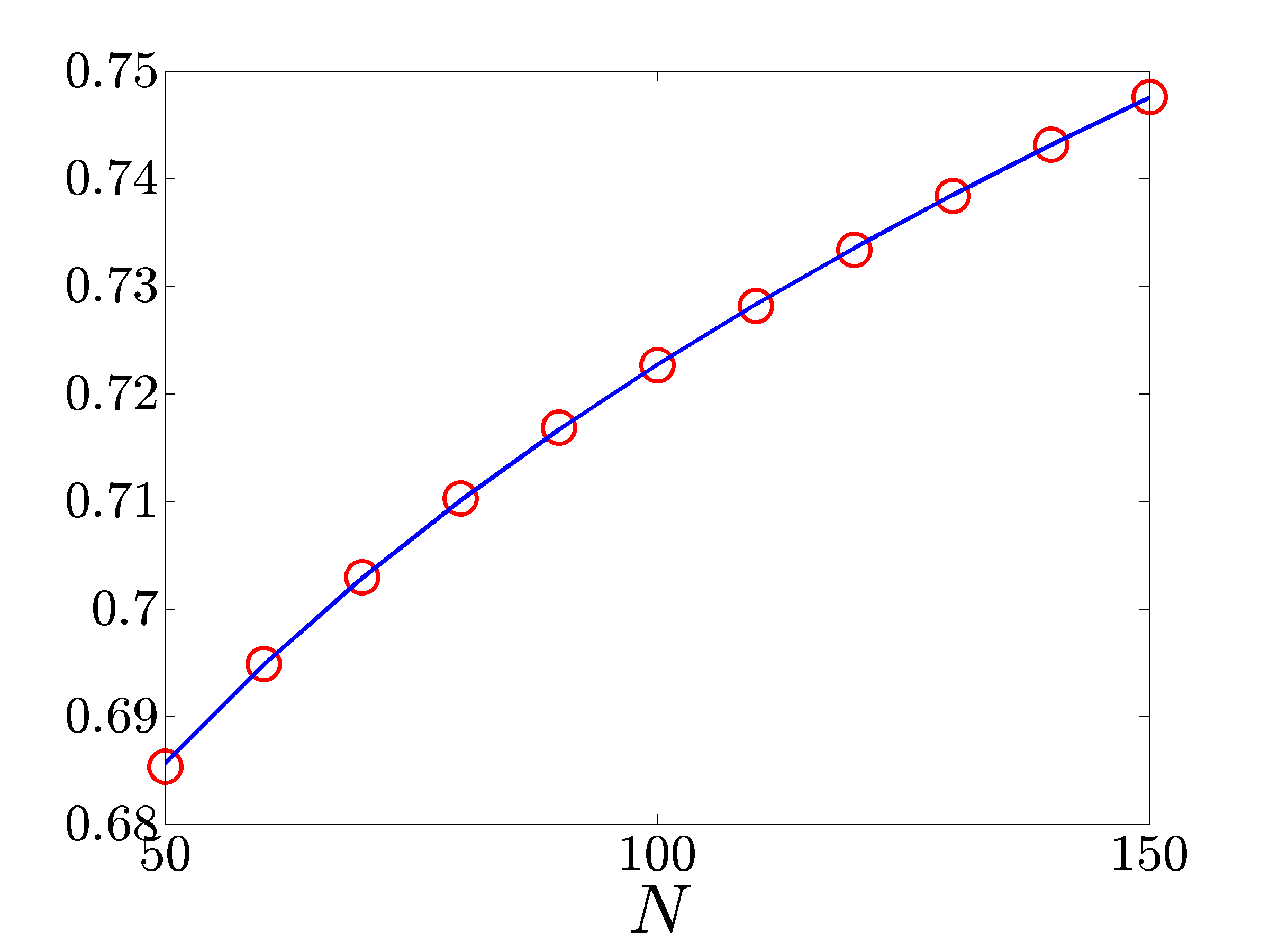}}
        \subfloat[$\Pi_l$]
        {\label{fig.Pi_l_Fg_2nd}
        \includegraphics[width=0.24\textwidth]{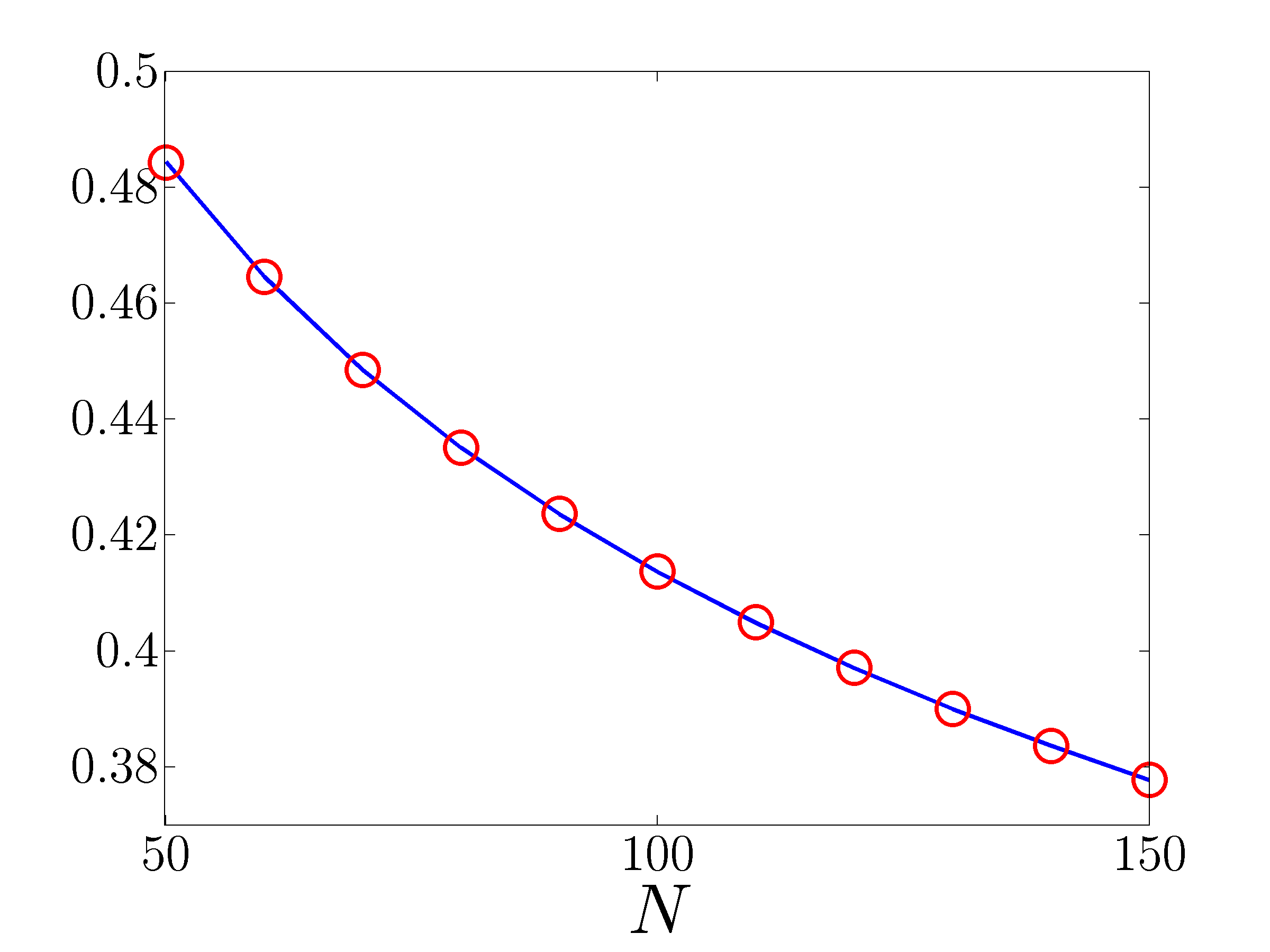}}
        \\
        \subfloat[$\Pctr$]
        {\label{fig.Pi_ctr_Fg_2nd}
        \includegraphics[width=0.24\textwidth]{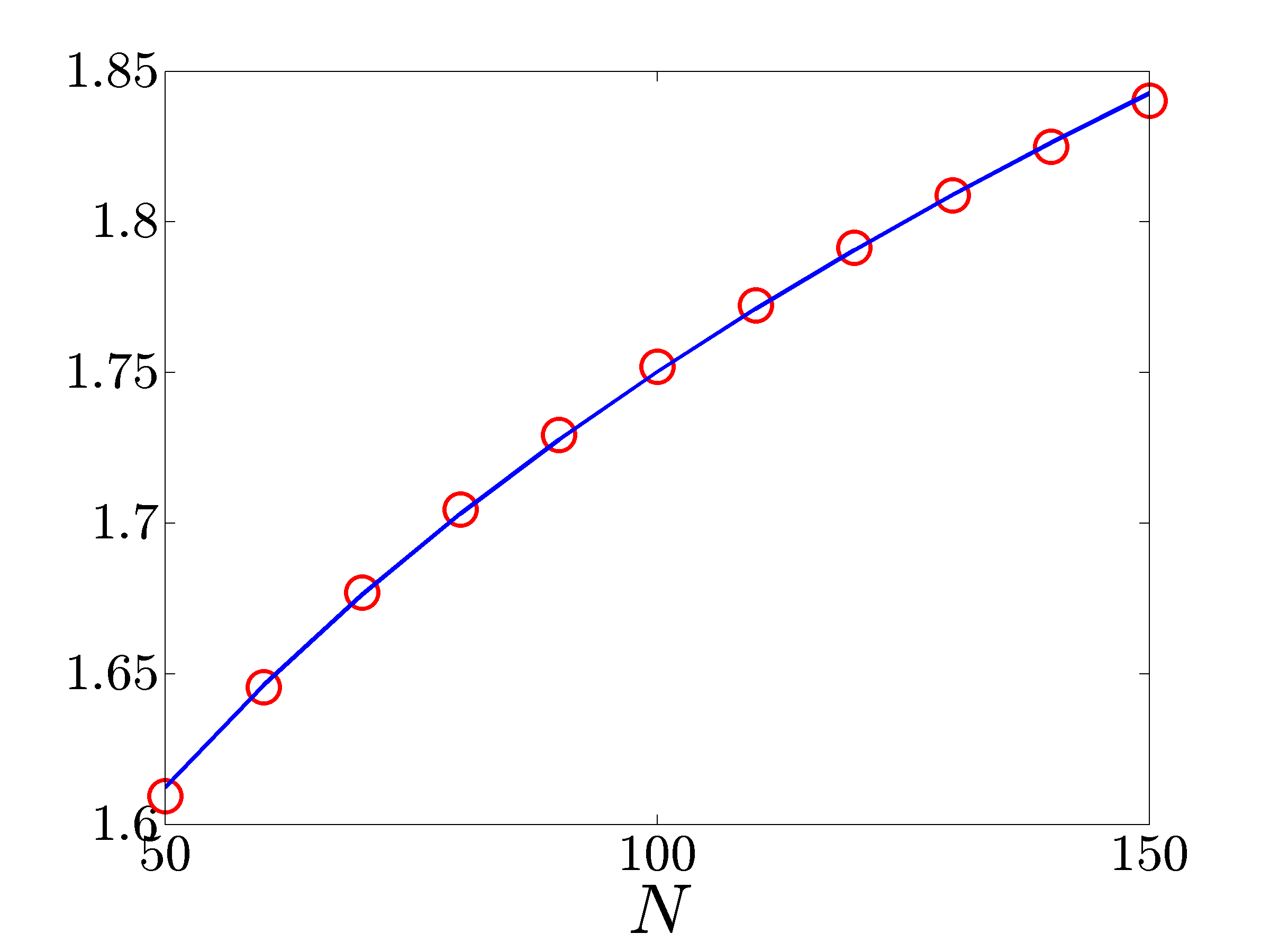}}
        \caption{Double-integrator model with the optimal non-symmetric gain
        obtained by solving~\ref{SH2} with $Q = I$ and $r = 1$ for formations with the fictitious follower:
        (a) $\Pi_g \, (\circ)$,         $0.0736 \sqrt[4]{N}+0.4900$ (curve)
        (b) $\Pi_l \, (\circ)$,         $1.1793 / \sqrt[4]{N} + 0.0408$ (curve);
        (c) $\Pctr \, (\circ)$, $0.2742 \sqrt[4]{N} + 0.8830$ (curve).}
    \label{fig.fourthrootfit_Fg_2nd}
    \end{figure}

As in the single-integrator model, our computations indicate that the optimal localized controller, obtained by solving~\ref{SH2} with $Q = I$ and $r = 1$, provides a fourth-root dependence of the macroscopic performance measure $\Pi_g$ on $N$; see Fig.~\ref{fig.Pi_g_Fg_2nd}. Furthermore, the microscopic performance measure and control energy asymptotically scale as $O(1/\sqrt[4]{N})$ and $O(\sqrt[4]{N})$, respectively; see Fig.~\ref{fig.Pi_l_Fg_2nd} and Fig.~\ref{fig.Pi_ctr_Fg_2nd}.

For comparison, we next provide the scaling trends of the performance measures for both the spatially uniform symmetric and look-ahead controllers. As in the single-integrator model, the spatially uniform symmetric gain
    $
    F_0
    =
    [ \, \alpha I \;\; \alpha I \;\; \beta I \, ]
    $
provides linear scaling of $\Pi_g$ with $N$ and the formation-size-independent $\Pi_l$ and $\Pctr$,
    \[
    \ba{rcl}
    \Pi_g(N)
    \!\! & = & \!\!
    ( N \,+\, 2 ) / ( 12 \alpha \beta )
    \,+\,
    1 / (2 \beta),
    \\[0.1cm]
    \Pi_l
    \!\! & = & \!\!
    1/(2 \alpha \beta)
    \,+\,
    1/(2 \beta),
    \\[0.1cm]
    \Pctr
    \!\! & = & \!\!
    \alpha / \beta
    \,+\,
    \beta / 2.
    \ea
    \]
On the other hand, for the double-integrator model the performance of the look-ahead strategy
    $
    K
    =
    F C
    =
    [\, \alpha C_f ~~  \beta I \,]
    $
heavily depends on the choices of $\alpha$ and $\beta$. In particular, for $\alpha = 1/4$ and $\beta = 1$, using similar techniques as in Section~\ref{sec.uni_nonsym_ctr}, we obtain
    \[
    \Pi_g(N)
    \,=\,
    \frac{1}{\sqrt{\pi}}
    \sum_{n\,=\,1}^N
    \dfrac{(N - n + 1)}{2 N \, \Gamma(2n)}
    \Big(
    8 \, \Gamma(2n - \dfrac{1}{2})
    \,+\,
    \Gamma(2n - \dfrac{3}{2})
    \Big),
    \]
which asymptotically leads to the formation-size-independent scaling of $\Pctr$ and the square-root scaling of $\Pi_g$ with $N$, i.e.,
    $
    \ds{\lim_{N \to \infty}}
    \Pi_g(N) / \sqrt{N}
    =
    16 / ( 3 \sqrt{2 \pi} )
    $.
This is in sharp contrast to $\alpha = \beta = 1$ which leads to an {\em exponential\/} dependence of $\Pi_g$ on $N$. Therefore, the design of the look-ahead strategy is much more subtle for double-integrators than for single-integrators.

\remark
For the double-integrator model with $K = F C = [\, K_p \;\; \beta I \,]$ and fixed $\beta > 0$ we next show convexity of~\ref{SH2} with respect to $K_p = K_p^T > 0$. The Lyapunov equation in~\ref{SH2}, for the block diagonal state weight $Q$ with components $Q_1$ and $Q_2$, can be rewritten in terms of the components of
    $
    P
    =
    \left[
    \ba{cc}
    P_1    &  P_0 \\
    P_0^T  &  P_2
    \ea
    \right],
    $
    \begin{subequations}
    \begin{align}
    \label{eq.P0}
      K_p P_0^T \,+\, P_0 K_p
      \,&=\,
      Q_1 \,+\, K_p K_p, \\
      \label{eq.P1}
      K_p P_2 \,-\, P_1 \,+\, \beta P_0
      \,&=\,
      \beta K_p, \\
      \label{eq.P2}
      2 \beta P_2
      \,&=\,
      P_0 \,+\, P_0^T \,+\, Q_2 \,+\, \beta^2 I.
    \end{align}
    \end{subequations}
Linearity of the trace operator in conjunction with $B_1 = [\, O ~~ I \,]^T$ and~(\ref{eq.P2}) yields
    \[
    \ba{rcl}
    J
    \!\! &=& \!\!
    \trace \left( P_2 \right)
    \,=\,
    \trace
    \left(
    2 P_0 + Q_2 + \beta^2 I
    \right)/(2\beta)
    \\[0.1cm]
    \!\! &=& \!\!
    \trace
    \left(
    K_p^{-1} Q_1 + K_p + Q_2 + \beta^2 I
    \right)/(2\beta),
    \ea
    \]
where the last equation is obtained by multiplying (\ref{eq.P0}) from the left
with $K^{-1}_p$ and using $\trace \, (K_p^{-1} P_0 K_p) = \trace \, (P_0)$. For $Q_1 \geq 0$, similar argument as in Section~\ref{sec.sym} can be used to conclude convexity of $J$ with respect to $K_p = K_p^T > 0$.

\section{Concluding remarks}
    \label{sec.con_rem}

We consider the optimal control of one-dimensional formations with nearest neighbor interactions between the vehicles. We formulate a structured optimal control problem in which local information exchange of relative positions between immediate neighbors imposes structural constraints on the feedback gains. We study the design problem for both the single- and the double-integrator models and employ a homotopy-based Newton's method to compute the optimal structured gains. We also show that design of symmetric gains for the single-integrator model is a convex optimization problem, which we solve analytically for formations with no fictitious followers. For double-integrators, we identify a class of convex problems by restricting the controller to symmetric position and uniform diagonal velocity gains. Furthermore, we investigate the performance of the optimal controllers by examining the asymptotic scalings of formation coherence and control energy with the number of vehicles.

For formations in which all vehicles have access to their own velocities, the optimal structured position gains for single- and double-integrators are similar to each other. Since these two models exhibit the same asymptotic scalings of global, local, and control performance measures, we conclude that the single-integrator model, which lends itself more easily to analysis and design, captures the essential features of the optimal localized design. We note that the tools developed in this paper can also be used to design optimal structured controllers for double-integrators with relative position and velocity measurements; this is a topic of our ongoing research.

As in~\cite{barmehhes09}, we employ perturbation analysis to determine the departure from a stabilizing spatially uniform profile that yields nominal diffusion dynamics on a one-dimensional lattice; in contrast to~\cite{barmehhes09}, we find the `mistuning' profile by optimizing a performance index rather than by performing spectral analysis. We also show how a homotopy-based Newton's method can be employed to obtain non-infinitesimal variation in feedback gains that minimizes the desired objective function. Furthermore, we establish several explicit scaling relationships and identify a spatially uniform non-symmetric controller that performs better than the optimal symmetric spatially varying controller ($O(\sqrt{N})$ vs.\ $O(N)$ scaling of coherence with $O (1)$ control energy per vehicle). This suggests that departure from symmetry can improve coherence of large-scale formations and that the controller structure may play a more important role than the optimal feedback gain design. On the other hand, our results demonstrate that the best performance is achieved with the optimal localized controller that is both non-symmetric and spatially-varying.

Currently, we are considering the structured feedback design for formations on general graphs~\cite{faxmur04,lafwilcauvee05,xiaboykim07,borkev08,tun09} with the objective of identifying topologies that lead to favorable system-theoretic properties~\cite{youscaleo10,zelmes11,linfarjovCDC10}. Even though this paper focuses on the optimal local feedback design for one-dimensional formations with path-graph topology, the developed methods can be applied to multi-agent problems with more general network topologies.

\section*{Acknowledgements}

The authors would like to thank anonymous reviewers and the associate editor for their valuable comments.

\appendix

\subsection{Gradient method for {\em \ref{SG}}}
    \label{app.grad_sym_K}

We next describe the gradient method for solving~\ref{SG}. Let us denote $k = [\,k_1 \, \cdots \, k_{N+1}\,]^T$. Starting with an initial guess $k^0$ that guarantees positive definiteness of $K^0$, vector $k$ is updated
    $
        k^{i+1}
        =
        k^i
        -
        s^i \, \nabla J(k^i),
    $
until the norm of gradient is small enough, $\| \nabla J(k^i) \| < \epsilon$. Here, $s^i$ is the step-size determined by the backtracking line search~\cite[Section~9.2]{boyvan04}: let $s^i=1$ and repeat $ s^i := \beta s^i$ with $\beta \in (0,1)$ until a sufficient decrease in the objective function is achieved,
    \[
       J( k^i - \, s^i \nabla J(k^i) )
       \, < \,
       J(k^i)
       \, - \,
       \alpha \,s^i\, \| \nabla J(k^i) \|^2 ,
    \]
where $\alpha \in (0,0.5)$. Note that $J(k^i)$ is defined as infinity if $K$ in~(\ref{eq.K}) determined by $k^{i}$ is not positive definite.
For $Q = I$,
    \begin{align*}
        J
        \,&=\,
        \dfrac{1}{2} \, \trace \left( K^{-1} + r K \right)
        \\
        \,&=\,
        \dfrac{1}{2} \sum_{n\,=\,1}^N
        \frac{\gamma_n(\gamma_{N+1}-\gamma_n)}{\gamma_{N+1}}
        \,+\,
        r
        \left(
        \dfrac{k_1 + k_{N+1}}{2}
        \,+\,
        \sum_{n\,=\,2}^N k_n
        \right),
    \end{align*}
the entries of the gradient $\nabla J$ are given by
    \beq
    \ba{rcl}
        \dfrac{\partial J}{\partial k_n}
        \!\! & = & \!\!
        r
         -
        \dfrac{1}{2}
        \,
        \ds{\sum_{i\,=\,1}^{n-1}}
        \Big( \frac{\gamma_i}{k_n \gamma_{N+1}} \Big)^2
         -
        \dfrac{1}{2}
        \,
        \ds{\sum_{i\,=\,n}^{N}}
        \Big(
        \frac{\gamma_{N+1} - \gamma_i}{k_n \gamma_{N+1}}
        \Big)^2,
        \\[0.45cm]
        & &
        n \in \{ 2,\ldots,N \},
        \\
        \dfrac{\partial J}{\partial k_1}
        \!\! & = & \!\!
        \dfrac{r}{2}
        \,-\,
        \dfrac{1}{2}
        \,
        \ds{\sum_{n\,=\,1}^N}
        \left( \frac{\gamma_{N+1}-\gamma_n}{k_1\gamma_{N+1}}
        \right)^2,
        \\
        \dfrac{\partial J}{\partial k_{N+1}}
        \!\! & = & \!\!
        \dfrac{r}{2}
        \,-\,
        \dfrac{1}{2}
        \,
        \ds{\sum_{n\,=\,1}^N}
        \left(
        \frac{\gamma_n}{k_{N+1}\gamma_{N+1}}
        \right)^2.
    \ea
    \non
    \eeq

\subsection{Performance of look-ahead strategy}
    \label{app.nonsymmetric}

We next derive the analytical expressions for the performance measures $\Pi_g$, $\Pi_l$, and $\Pctr$ obtained with the look-ahead strategy for the single-integrator model. The solution of the Lyapunov equation~(\ref{eq.Lctrb}) with $F C = \alpha C_f$ is determined by~(\ref{eq.Lint}). Since the $i$th entry of the first column of the lower triangular Toeplitz matrix \mbox{$(sI + \alpha C_f)^{-1}$} is $\alpha^i / (s + \alpha)^i$, the corresponding entry of the matrix exponential in~(\ref{eq.Lint}) is determined by the inverse Laplace transform of $\alpha^i / (s + \alpha)^i$,
    \[
    \alpha
    \left( \alpha t \right)^{i-1}
    \mre^{- \alpha t}/(i-1)!.
    \]
Thus, the $n$th element on the main diagonal of the matrix $L$ in~(\ref{eq.Lint}) is given by
    \beq
    \label{eq.Lnn}
    \ba{rcl}
    L_{nn}
    \!\! & = & \!\!
    \ds{\int_{0}^{\infty}}
    \ds{\sum_{i\,=\,1}^{n}}
    \left(  \alpha \mre^{- \alpha t} \, \dfrac{( \alpha t)^{i-1}}{(i-1)!} \right)^2
    \mrd t
    \\[0.45cm]
    \!\! & = & \!\!
    \dfrac{\alpha \, \Gamma(n + 1/2)}{\sqrt{\pi} \, \Gamma(n)}
    \,=\,
    \dfrac{ \alpha \, (2n)!}{2^{2n} (n-1)! \, n!},
    \ea
    \eeq
thereby yielding
    \beq
    \label{eq.pig}
    \Pi_g
    =
    \sum_{n\,=\,1}^N
    \dfrac{L_{nn}}{N}
    =
    \dfrac{2 \, \alpha \, \Gamma(N + 3/2)}{3 \, \sqrt{\pi} \, \Gamma(N + 1)}
    =
    \frac{2}{3}
    \frac{ \alpha \, (2N+2)! }{ 2^{2N+2} N! (N+1)!}.
    \eeq
A similar procedure can be used to show that the $n(n+1)$th entry of $L$ is determined
    \beq
    \label{eq.Lnnp1}
    L_{n(n+1)}
    \,=\,
    L_{(n+1)(n+1)}
    \,-\,
     \alpha /2
    ,
    ~~
    n \,=\, 1,\ldots,N-1.
    \eeq
Now, from~(\ref{eq.Lnnp1}) and the fact that $L_{11} =  \alpha /2$ we obtain
    \[
    \Pi_l
    =
    \frac{1}{N}
    \,
    \trace
    \left(
    T L
    \right)
    =
    \frac{2}{N}
    \left(
    \sum_{n\,=\,1}^N  L_{nn}
    \,- \,
    \sum_{n\,=\,1}^{N-1} L_{n(n+1)}
    \right)
    =
    \alpha,
    \]
Similarly,
    \[
    \ba{rcl}
    \Pctr
    \!\! & = & \!\!
    (1/N)
    \,
    \trace
    \left(
    L C_f^T C_f
    \right)
    \\[0.25cm]
    \!\! & = & \!\!
    \dfrac{2}{N}
    \left(
    \ds{\sum_{n\,=\,1}^N}  L_{nn}
    \,- \,
    \ds{\sum_{n\,=\,1}^{N-1}} L_{n(n+1)}
    \right)
    \, - \,
    \dfrac{1}{N} \, L_{NN}
    \\[0.45cm]
    \!\! & = & \!\!
    \alpha
    \,-\,
    (1/N) \, L_{NN}.
    \ea
    \]
Using Stirling's approximation
    $
    n!
    \approx
    \sqrt{2\pi n}
    \,
    ( n / \mre )^n
    $
for large $n$, we have
    \begin{align*}
    \lim_{n \, \to \, \infty}
    \frac{L_{nn}}{\sqrt{n}}
    \, = \,
    \lim_{n \, \to \, \infty}
    \frac{ \alpha }{\sqrt{\pi}}
    \,
    \sqrt{\frac{ n}{ n-1 }}
    \left(
    \frac{n}{n-1}
    \right)^{n-1}
    \frac{1}{\mre}
    \,=\,
    \frac{ \alpha }{\sqrt{\pi}},
    \end{align*}
where we used the fact that
    $
    \ds{\lim_{n \, \to \, \infty}}
    \left(
    n/(n-1)
    \right)^{n-1}
    \,=\,
    \mre.
    $
Consequently,
    $
    \ds{\lim_{N \,\to\,\infty}}
    \Pctr(N)
    =
    \alpha.
    $
From (\ref{eq.Lnn}) and (\ref{eq.pig}), it follows that $\Pi_g = (2/3) L_{(N+1)(N+1)}$ and thus,
    $
    \ds{\lim_{N \,\to\,\infty}}
    \Pi_g(N) / \sqrt{N}
    =
    (2 \alpha) / (3 \sqrt{\pi}).
    $
We conclude that $\Pi_g$ asymptotically scales as a square-root function of $N$ and that $\Pctr$ is formation-size-independent as $N$ increases to infinity.


    \vspace*{-1.5cm}
\begin{biography}[{\includegraphics[width=0.85in,height =1.1in,clip,keepaspectratio]{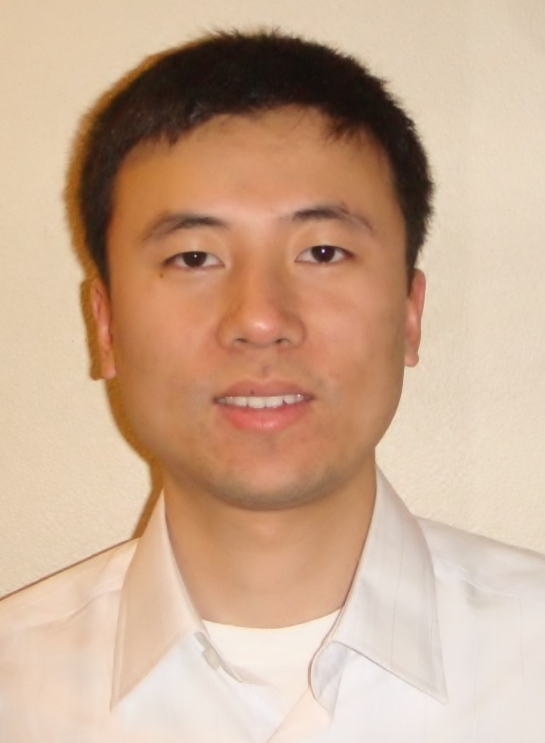}}]{Fu Lin} (S'06) received his Bachelor of Science degree in Instrument Science
and Engineering from Shanghai Jiaotong University in 2005. Currently, he is a Ph.D.\
candidate in the Department of Electrical and Computer Engineering at the University of Minnesota,
Minneapolis. His primary research interests are in the analysis and design of
optimal distributed controllers using tools from convex optimization, compressive sensing, and graph theory.
\end{biography}
     \vspace*{-1.95cm}
\begin{biography}[{\includegraphics[width=1in,height =1.25in,clip,keepaspectratio]{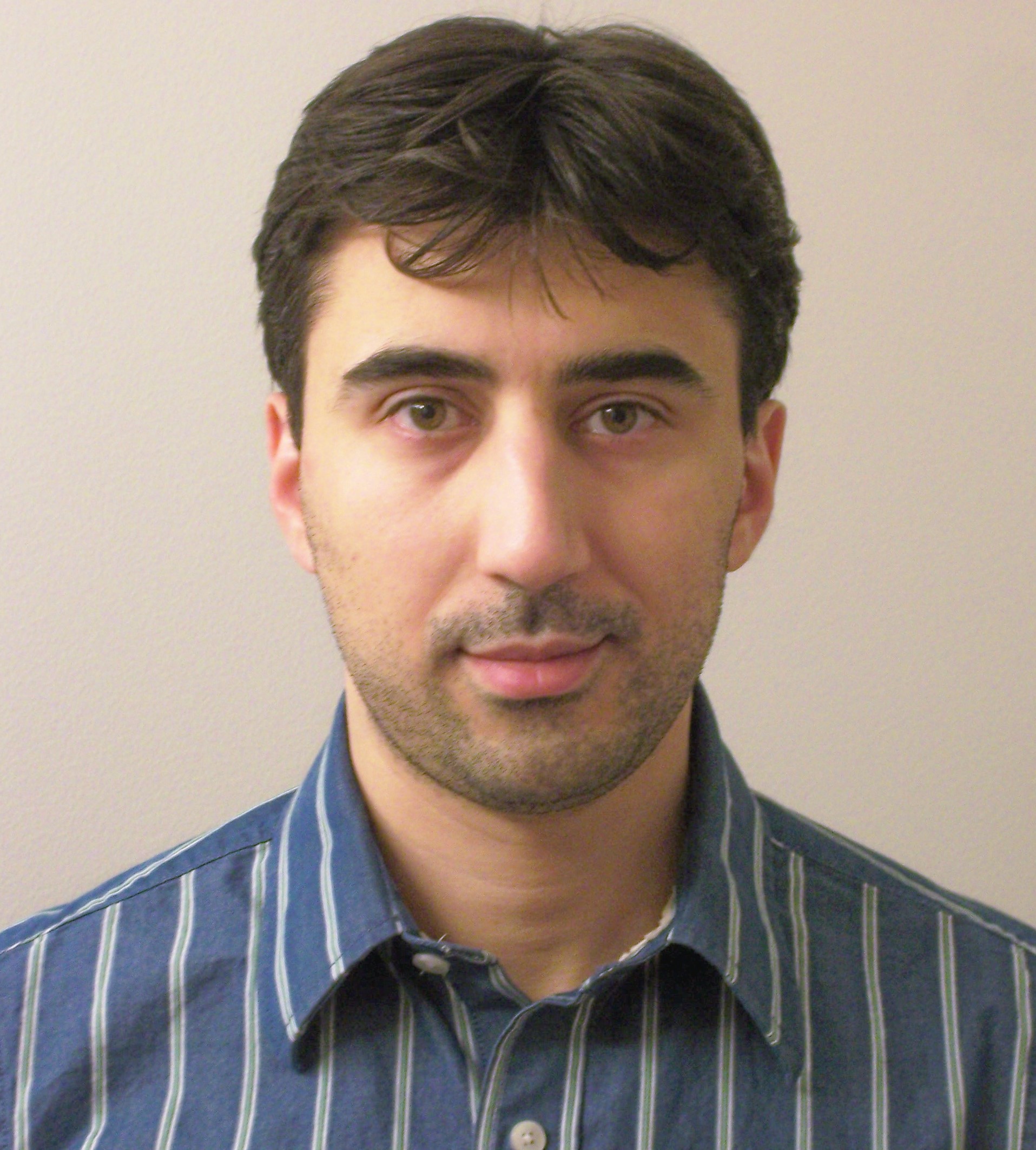}}]{Makan Fardad} received the B.S.\ and M.S.\ degrees in electrical engineering from Sharif University of Technology and Iran University of Science and Technology, respectively. He received the Ph.D.\ degree in mechanical engineering from the University of California, Santa Barbara, in 2006. He was a Postdoctoral Associate at the University of Minnesota, Minneapolis, before joining the Department of Electrical Engineering and Computer Science at Syracuse University as an Assistant Professor in August 2008. His research interests are in modeling, analysis, and optimal control of distributed and large-scale interconnected systems.
\end{biography}
     \vspace*{-1.95cm}
\begin{biography}[{\includegraphics[width=1in,height =1.25in,clip,keepaspectratio]{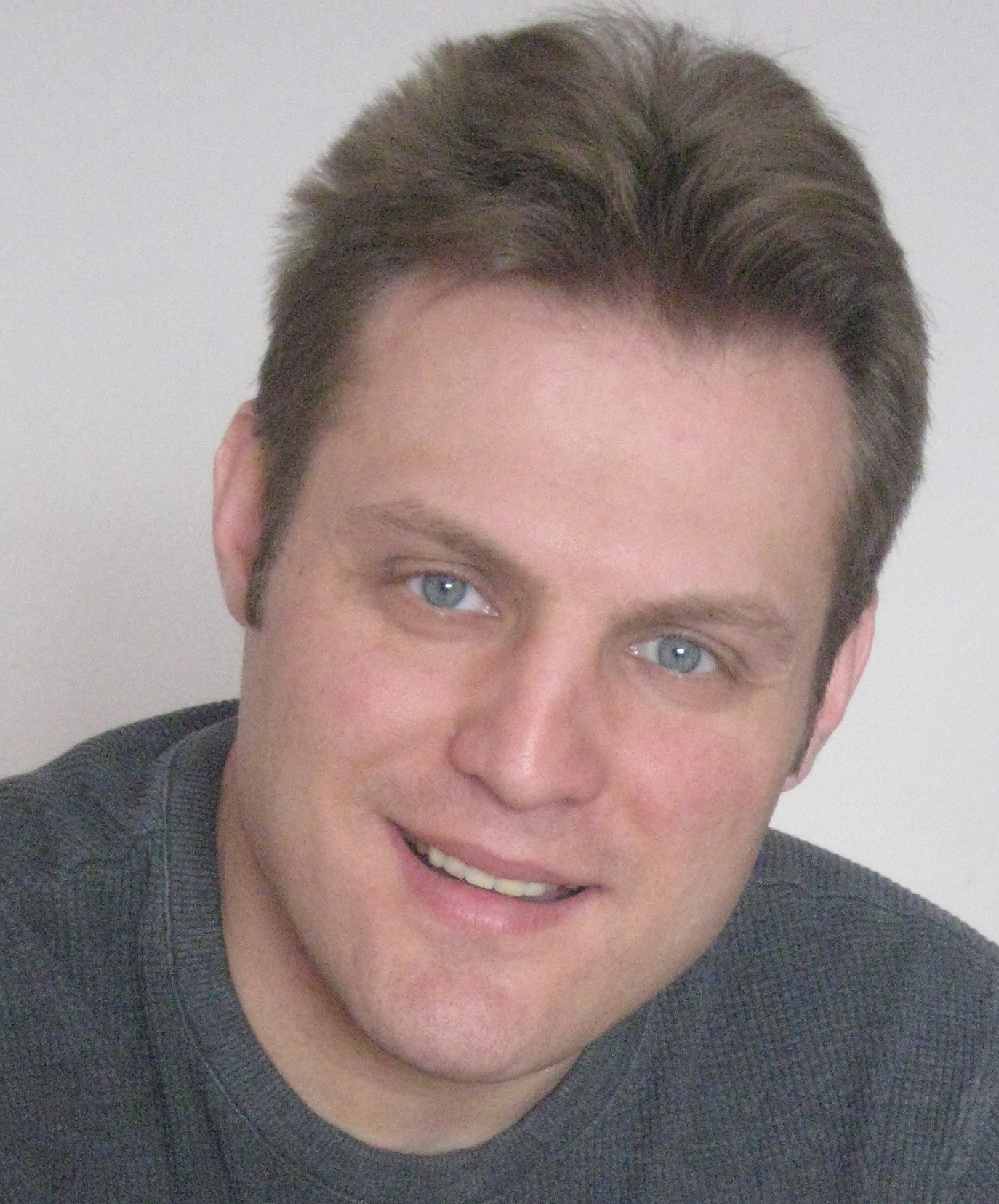}}]{Mihailo R.\ Jovanovi\'c}  (S'00--M'05)
received the Dipl.\ Ing.\ and M.S.\ degrees from the University of Belgrade, Serbia, in 1995 and 1998, respectively, and the Ph.D.\ degree from the University of California, Santa Barbara, in 2004. Before joining the University of Minnesota, Minneapolis, he was a Visiting Researcher with the Department of Mechanics, the Royal Institute of Technology, Stockholm, Sweden, from September to December 2004. Currently, he is an Associate Professor of Electrical and Computer Engineering at the University of Minnesota, where he serves as the Director of Graduate Studies in the interdisciplinary Ph.D.\ program in Control Science and Dynamical Systems.

Dr.\ Jovanovi\'c's expertise is in modeling, dynamics, and control of large-scale and distributed systems and his current research focuses on sparsity-promoting optimal control, dynamics and control of fluid flows, and fundamental limitations in the control of vehicular formations. He is a member of APS and SIAM and has served as an Associate Editor of the IEEE Control Systems Society Conference Editorial Board from July 2006 until December 2010. He received a CAREER Award from the National Science Foundation in 2007, and an Early Career Award from the University of Minnesota Initiative for Renewable Energy and the Environment in 2010.

\end{biography}

\end{document}

%% file: commands.tex
\newtheorem{remark}{Remark}

\newcommand{\enma}[1]   {\ensuremath{#1}}

\newcommand{\nn}{\nonumber}

\newcommand{\beq}{\begin{equation}}
\newcommand{\eeq}{\end{equation}}
\newcommand{\bseq}{\begin{subequations}}
\newcommand{\eseq}{\end{subequations}}
\newcommand{\beqn}{\begin{eqnarray}}
\newcommand{\eeqn}{\end{eqnarray}}
\newcommand{\ba}{\begin{array}}
\newcommand{\ea}{\end{array}}
\newcommand{\bct}{\begin{center}}
\newcommand{\ect}{\end{center}}
\newcommand{\btmz}{\begin{itemize}}
\newcommand{\etmz}{\end{itemize}}
\newcommand{\benum}{\begin{enumerate}}
\newcommand{\eenum}{\end{enumerate}}

\newcommand{\diag}      {\enma{\mathrm{diag}}}

\newcommand{\trace}     {\enma{\mathrm{trace}}}

\newcommand{\matbegin}{
        \left[
}
\newcommand{\matend}{
        \right]
}

\newcommand{\tbo}[2]{
  \matbegin \begin{array}{c}
       #1 \\ #2
       \end{array} \matend }

\newcommand{\obt}[2]{
  \matbegin \begin{array}{cc}
       #1 & #2
       \end{array} \matend }
\newcommand{\obth}[3]{
  \matbegin \begin{array}{ccc}
       #1 & #2 & #3
       \end{array} \matend }
\newcommand{\tbt}[4]{
  \matbegin \begin{array}{cc}
       #1 & #2 \\ #3 & #4
       \end{array} \matend }

\newcommand{\be}{\begin{equation}}
\newcommand{\ee}{\end{equation}}

\newcommand{\cplxs}{ C\kern -.35em \rule{0.03 em}{.7 ex}~   }

\def\complex{\hbox{C\kern -.45em \rule{0.03 em}{1.5 ex}}~}

\newcommand{\bi}{\begin{itemize}}
\newcommand{\ei}{\end{itemize}}